\documentclass[12pt,reqno]{amsart}
\usepackage{amsmath,amsthm,amsbsy,amssymb, longtable,hhline}
\topmargin by-\headheight \advance
\topmargin by-\headsep
\evensidemargin
\oddsidemargin \marginparwidth 0.5in
\usepackage{amsmath, amssymb, graphicx, array, verbatim, psfig,
caption, ccaption, epsfig, xspace}

\newtheorem{theorem}{Theorem}
\newtheorem{lemma}{Lemma}
\newtheorem{proposition}{Proposition}
\newtheorem{corollary}{Corollary}

\usepackage{footnpag}

\newcommand{\Area}{\mathrm{Area}}
\newcommand{\length}{\mathrm{length}}

\newcolumntype{C}{>{$}c<{$}}
\newcolumntype{L}{>{$}l<{$}}
\newcolumntype{R}{>{$}r<{$}}

\def\A{\mathcal A} 

\newcommand{\C}{{\mathcal{C}}}

\newcommand{\T}{\mathcal{T}}
\newcommand{\B}{\mathcal{B}}
\newcommand{\D}{\mathcal{D}}
\newcommand{\E}{\mathcal{E}}

\newcommand{\I}{\mathcal{I}}
\newcommand{\LL}{\mathcal{L}}
\newcommand{\U}{\mathcal{U}}

\newcommand{\M}{\mathcal{M}}
\renewcommand{\P}{\mathcal{P}}

\newcommand{\kk}{{\mathbf k}}
\def\g{\gamma}

\def\RR{\mathbb{R}}
\def\NN{\mathbb{N}}

\def\ZZ{\mathbb{Z}}
\newcommand{\F}{{\mathfrak{F}}}
\newcommand{\FQ}{{\mathfrak{F}_{_Q}}}
\newcommand{\FQI}{{\mathfrak{F}^\I_{_Q}}}

\newcommand{\FQo}{{\mathfrak{F}_{_{Q,\mathrm{odd}}}}}
\newcommand{\FQoI}{{\mathfrak{F}^\I_{_{Q,\mathrm{odd}}}}}
\newcommand{\FQe}{{\mathfrak{F}_{_{Q,\mathrm{even}}}}}
\newcommand{\FQeI}{{\mathfrak{F}^\I_{_{Q,\mathrm{even}}}}}

\newcommand{\TT}{{\mathtt{T}}}

\renewcommand{\k}{{\mathbf{k}}}

\newcommand{\Tk}{{\mathcal{T}_{_{\k}}}}

\newcommand{\Noo}{ N_{\mathrm{odd},\mathrm{odd}}}
\newcommand{\Noe}{ N_{\mathrm{odd},\mathrm{even}}}

\newcommand{\Neo}{ N_{\mathrm{even},\mathrm{odd}}}

\newcommand{\huu}{h_1^{\mathrm u}}
\newcommand{\hdu}{h_2^{\mathrm u}}
\newcommand{\htu}{h_3^{\mathrm u}}
\newcommand{\hud}{h_1^{\mathrm d}}
\newcommand{\hdd}{h_2^{\mathrm d}}
\newcommand{\htd}{h_3^{\mathrm d}}

\newcommand{\vfi}{\varphi}
\newcommand{\vfit}{\widetilde\varphi}
\def\z{\overline z}
\def\W{\mathfrak G}
\def\wU{\widetilde\U}

\begin{document}

\author[C. Cobeli and A. Zaharescu]
{Cristian Cobeli and Alexandru Zaharescu}\footnotetext{CC is partially supported by the CERES Programme of the Romanian Ministry of Education and Research, contract 4-147/2004.}

\address{
CC:
Institute of Mathematics of the Romanian Academy,
P.O. Box \mbox{1-764}, Bucharest 70700,
Romania.}
\email{cristian.cobeli@imar.ro}

\address{
AZ:
Institute of Mathematics of the Romanian Academy,
P.O. Box \mbox{1-764}, Bucharest 70700,
Romania.}

\curraddr{
AZ:
Department of Mathematics,
University of Illinois at Urbana-Champaign,
Altgeld Hall, 1409 W. Green Street,
Urbana, IL, 61801, USA.}
\email{zaharesc@math.uiuc.edu}

\title[A density theorem on even Farey fractions]
{A density theorem on even Farey fractions}


\begin{abstract}
Let $\FQ$ be the Farey sequence of order $Q$ 
and let $\FQo$ and $\FQe$ be the set of those Farey fractions of order $Q$
with odd, respectively even denominators.
A fundamental property of $\FQ$ says that the sum of
denominators of any pair of neighbor fractions is always greater than
$Q$. This property fails for $\FQo$ and for $\FQe$.
The local density, as $Q\to\infty$, of the normalized pairs 
$(q'/Q,q''/Q)$, where $q',q''$ are denominators of consecutive 
fractions in $\FQo$, was computed in ~\cite{Integers}. The density
increases over a series of quadrilateral
steps ascending in a harmonic series towards the point
$(1,1)$. Numerical computations for small values of $Q$ 
suggest that such a result should rather
occur in the even case, while in the odd case the distribution of
the corresponding points appears to be more uniform. Reconciling with 
the numerical experiments, in this paper we show that, as $Q\to\infty$, 
the local densities in the odd and even case coincide.
\end{abstract}

\maketitle

\section{Introduction}\label{Section1}

Questions concerned with Farey sequences have a long history.
In some problems, such as for instance those related to the 
connection between Farey fractions and Dirichlet $L-$functions, 
one is lead to consider subsequences of Farey fractions defined 
by congruence constraints. 
Recently it has been realized that
knowledge of the distribution
of subsets of Farey fractions with congruence constraints
would also be useful in the study of the
periodic two-dimensional Lorentz gas. This is a billiard
system on the two-dimensional torus with one or more
circular regions (scatterers) removed 
(see~\cite{Si},~\cite{Bu},~\cite{Ch},~\cite{BZ}). 
Such systems were introduced in 1905 by 
Lorentz~\cite{L} to describe the dynamics of
electrons in metals. A problem raised by Sina\u {\i}
on the distribution of the free path length for this
billiard system, when small scatterers are placed at integer
points and the trajectory of the particle
starts at the origin, was solved in~\cite{BGZ1}, \cite{BGZ2},
using techniques developed in~\cite{ABCZ}, \cite{BCZ1}, 
\cite{BCZ2} to study the local spacing distribution of Farey 
sequences. 

The more general case when the trajectory starts 
at a given point with rational coordinates is intrinsically 
connected with the problem of the distribution of
Farey fractions satisfying congruence constraints. 
For example, the case when the trajectory starts from the center 
$(1/2,1/2)$ of the unit square
is related to the distribution of Farey fractions
with odd numerators and denominators.
The distribution of the free path length computed in~\cite{BGZ1} 
and~\cite{BGZ2} is totally different from
the one obtained in~\cite{BZ}. This confirms the
intuition of physicists that, unlike in the case when the
trajectory starts at the origin, if one averages over
the initial position of the particle, the distribution
will have a tail. It is then reasonable to expect
that, in terms of the distribution of Farey fractions,
new phenomena would be encountered when one replaces
the entire sequence of Farey fractions by a subsequence
defined by congruence constraints. 

As we shall see below,
already the case of the subsequence of Farey fractions
with even denominators presents nontrivial complications.
This is mainly due to the fact that in $\FQ$ there is a 
large number of tuples of consecutive fractions with odd denominators
and length growing to infinity with $Q$. 
Recently some questions on the distribution of Farey fractions with
odd denominators have been investigated in~\cite{BCZ3}, 
\cite{Integers}, and~\cite{H}. It is the purpose of
this work to derive a result on fractions with even denominators.

Two fundamental properties of the Farey sequence of order $Q$
state that if $a'/q'<a''/q''$ are consecutive
elements of $\FQ$, then $a''q'-a'q''=1$, and $q'+q''>Q$. 
These properties play an
essential role in questions concerned with the distribution
of Farey fractions. In fact, in any problem where one has
an element $a'/q'$ of $\FQ$ and needs to find the next
element of $\FQ$, call it $a''/q''$, one can use the
above two properties in order to determine $a''/q''$,
as follows. The equality $a''q'-a'q''=1$ uniquely determines 
$a''$ in terms of $a', q'$ and $q''$. In order to find
$q''$ in terms of $a'$ and $q'$, notice from the above 
equality that $a'q''\equiv -1$ (mod $q'$).
The inequalities $q'+q''>Q$ and $q''\le Q$ show that $q''$ belongs 
to the interval $(Q-q',Q]$, which
contains exactly one integer from each residue class modulo $q'$.
Only one of these integers satisfies the congruence
$a'q''=-1$ (mod $q'$), and this uniquely determines $q''$
in terms of $q'$ and $a'$.
Complications arise when one studies a subsequence of Farey 
fractions with denominators in a given residue class modulo 
an integer number $d\ge 2$, since in such a case the above
two properties fail (see~\cite{BCZ3}, \cite{Integers}, \cite{H} 
for the case of fractions with odd denominators). 

\noindent
\begin{figure}[t]
    \begin{minipage}[t]{0.48\linewidth}
      \centering
        \includegraphics*[width=0.98\linewidth]{desen.12}
        \caption{\small From light to dark are represented the sets of
                type $\TT(1)$.}\label{Figure2}
     \end{minipage}
     \hfill
    \begin{minipage}[t]{0.48\textwidth}
      \centering
        \includegraphics*[width=0.98\linewidth]{desen.13}
                \caption{\small The sets of type $\TT(3)$, $\TT(4)$ and $\TT(r)$, 
                                with $r\ge 5$.}\label{Figure3}
    \end{minipage}
\end{figure}
\indent
In the present paper we
study the relative size of consecutive even denominators
in Farey series. Although the inequality $q'+q''>Q$
fails in this case too, we shall see that the points
$(q'/Q,q''/Q)$ have a limiting distribution inside the
unit square $[0,1]\times[0,1]$, as
$Q\rightarrow\infty$.
In the following we let
  \begin{equation*}
    \begin{split}
     \FQe=\Big\{\frac aq\colon\ 1\le a\le q\le Q, \ 
           \gcd(a,q) = 1,\ q\equiv 0\pmod 2\Big\}\,
    \end{split}
  \end{equation*}
and we always assume that the elements of $\FQe$ are arranged in 
increasing order. We call a Farey fraction {\em odd} if its denominator 
is odd and {\em even} if its denominator is even, respectively.
A new feature that the sequence $\FQe$ brings in this type of 
problems comes from the following phenomenon. We know from the
equality $a''q'-a'q''=1$ that the denominators $q',q''$ of any two
consecutive fractions $a'/q'<a''/q''$ in $\FQ$ are relatively prime,
and in particular not both of them are even. Therefore, if we
study the subset $\FQo$ of odd Farey fractions in $\FQ$, we
know that any two consecutive elements of $\FQo$ are either
consecutive in $\FQ$, or there is exactly one element of
$\FQe$ between them. Thus there are only two types of situations
to consider. By contrast, we may have any number of elements
from $\FQo$ between two consecutive elements of $\FQe$.
This more complicated context that arises in the even case,
treated in the present paper, forces us to go through a
significantly larger amount of data than in the odd case.
This is also reflected in the larger variety of situations
that appear in
Figures~\ref{Figure2}-\ref{Figure5} and Tables~\ref{Table1}-\ref{Table2} 
below, which show various aspects of the distribution in this case.
In what follows we study the local density of points $(q'/Q,q''/Q)$, 
with $q',q''$ denominators of consecutive Farey fractions in $\FQe$,
which lie around a given point $(u,v)$ in the unit square.
We shall show that this local density approaches a certain limit 
$g(u,v)$ as $Q\rightarrow\infty$, and we provide 
an explicit formula for $g(u,v)$. 

\noindent
\begin{figure}[t]
    \begin{minipage}[t]{0.48\linewidth}
      \centering
        \includegraphics*[width=0.98\linewidth]{desen.14}
        \caption{\small The sequences of sets of type $\TT(2)$.}\label{Figure4}
     \end{minipage}
     \hfill
    \begin{minipage}[t]{0.48\textwidth}
      \centering
        \includegraphics*[width=0.98\linewidth]{desen.15}
        \caption{\small The covering of $\D(0,2)$ by the sets of all types.}\label{Figure5}
    \end{minipage}
\end{figure}
\indent
The property of Farey fractions that first drew attention two centuries
ago was their very uniform distribution of in $[0,1]$ (see the survey
paper~\cite{Survey} and the references therein). It is natural to
expect various distribution results, in particular the one obtained
in the present paper, to continue to hold in subintervals of $[0,1]$. 
Let $\I\subseteq [0,1]$ be a fixed subinterval, and denote 
$\FQI=\FQ\cap \I$ and $\FQeI=\FQe\cap \I$.
We let
  \begin{equation*}
    \begin{split}
       \D^\I_{Q,even}:=\big\{(q',q'')
\colon\ q',q''\ \text{denominators of consecutive
           fractions in}  \ \FQeI \big\}\,,
    \end{split}
  \end{equation*}
and consider the set $\D^\I_{Q,even}/Q$, which is a subset of
the unit square $[0,1]\times[0,1]$. 
If $\I$ is the complete interval $[0,1]$, we drop the superscript and write
$\D_{Q,even}=\D_{Q,even}^{[0,1]}$.

For each point 
$(u,v)\in[0,1]\times[0,1]$, we take a small square $\square$ centered
at $(u,v)$ and count the number of points from $\D^\I_{Q,even}/Q$
which lie inside the square $\square$. We shall see that, as
$Q\rightarrow\infty$, the proportion of points from $\D^\I_{Q,even}/Q$
that fall inside $\square$ approaches a certain limit. This limit
will be proportional to the area of the square $\square$. After 
dividing this limit by $\Area(\square)$ and letting the side of the 
square $\square$ tend to $0$, we arrive at a limit, call it
$g^\I(u,v)$, which will only depend on the point $(u,v)$, and possibly
on the interval $\I$. Thus, we let
    \begin{equation}\label{eqdefg}
      g^\I(u,v):=\lim_{\Area(\square)\to 0}
           \frac{\lim\limits_{Q\to \infty}
             \frac{\#\big(\square\cap\D^\I_{Q,even}/Q\big)}{\#\D^\I_{Q,even}}}
                  {\Area(\square)}\,,
    \end{equation}
in which $\square\subset\RR^2$ are squares centered at $(u,v)$.
We put $g(u,v)=g^{[0,1]}(u,v)$.
Theorem~\ref{Theorem1} below shows that the limiting local density
function $g(u,v)$ exists, and its value is calculated explicitly.

Since the Farey fractions are distributed in $[0,1]$ symmetrically
with respect to $1/2$, the components of the argument of the density
$g(u,v)$ will play a symmetric%
\footnote{We use tho word {\em symmetric} for tuples with components
listed in reverse order of one another, and also for points situated
symmetrically with respect to the first diagonal.}
 role and we shall have $g(u,v)=g(v,u)$, for any 
$(u,v)\in[0,1]^2$.
For convenience, in the following we shall use the symbol $z$ for either of
the variables $u$ or $v$ and $\z$ for the other. Also, to write
shortly the characteristic function of a system of {\em conditions} 
(equalities or inequalities in variables $u$ and $v$), 
we denote:
    \begin{equation*}
      \varphi(\text{conditions})=
        \begin{cases}
          1, & \text{if conditions hold true for 
                        $(z,\z)=(u,v)$ or $(z,\z)=(v,u)$};\\
          0, & \text{else,}
        \end{cases}
    \end{equation*}
and
    \begin{equation*}
      \vfit(\text{conditions})=
        \begin{cases}
          1, & \text{if conditions hold true for 
                        $(z,\z)=(u,v)$ and $(z,\z)=(v,u)$};\\
          0, & \text{else.}
        \end{cases}
    \end{equation*}

\begin{theorem}\label{Theorem1}
\noindent
The local density in the unit square of points 
$\left(q'/Q,q''/Q\right)$, where $q'$ and $q''$ are denominators of
neighbor fractions in $\FQe$, approaches a limiting density
$g$ as $Q\rightarrow\infty$. Moreover, for any
real numbers $u,v$ with $0\le u,v\le 1$, 
     \begin{equation}\label{eqT1} 
        \begin{split}
    g(u,v)=&\phantom{+}
                \sum_{j=1}^\infty 
        \frac 1{j}
                \vfit\Big(z<1;\ j<\frac{z+\z}{1-z}\Big)\\
                &+              \sum_{j=1}^\infty 
        \frac 1{2j}\bigg\{\vfi\Big((j+1)z+\z=j, \mathrm{\ if\ }
                \frac{j-1}{j+1}<z<\frac{j}{j+2}\Big)\\
        &\phantom{\sum_{\substack{k\ge 4\\k\text{ even}}}\frac{1}{2k}
                k=\frac{x+y}{1-z}\qquad\qquad}
        +\vfi\Big(z=1,\ \mathrm{\ if\ }
                \frac{j-1}{j+1}<\z<1\Big)\bigg\}\\
        &+              \sum_{j=1}^\infty  
         \bigg\{\frac{2j+1}{8j(j+1)}
                \vfi\Big(z=\frac{j-1}{j+1};\ \z=1\Big)\\
        &\phantom{\sum_{\substack{k\ge 4\\k\text{ even}}}\frac{1}{2k}++}
                +\frac{j+2}{4j(j+1)}
                \vfit\Big(z=\frac{j}{j+2}\Big)
                +\frac{1}{4j}
                \vfit\Big(z=1\Big)
        \bigg\}\,.
        \end{split}
     \end{equation}
\end{theorem}


Figures~\ref{Figure2}-\ref{Figure5} show how the unit square 
is covered by countably many polygons, on the interior of which
the local density function $g(u,v)$ is constant.
The explicit value of that constant is provided by
Theorem~\ref{Theorem1}. We remark that for a point $(u,v)$, which
is an interior point of one of the polygons that form
the above covering of the unit square, the second sum
and the third sum on the right side of \eqref{eqT1}
vanish, and in the first sum only finitely many terms
are nonzero, namely those corresponding to the values
of $j\geq1$ for which one has simultaneously $j<(u+v)/(1-v)$
and $j<(u+v)/(1-u)$. So, $g(u,v)$ reduces in this case to
a partial sum of the harmonic series, with more and more
terms of the series to be counted as the point $(u,v)$
is chosen closer and closer to the upper-right corner
$(1,1)$ of the unit square.
Additional terms, contained in the second and in the third
sum on the right side of \eqref{eqT1}, only appear in the case
when $(u,v)$ lies on one of the sides or coincides with
one of the vertices of one of the polygons that form
the covering of the unit square.

Turning to the same problem on shorter intervals, if $\I$
does not have $1/2$ as midpoint, than the symmetry
$g^\I(u,v)=g^\I(v,u)$ is apriori not at all obvious.
The next theorem provides
the stronger result that $g^\I(u,v)$ not only exists, for any
$\I\subseteq [0,1]$, but that it is independent of $\I$.

\begin{theorem}\label{TheoremI}
Let  $\I\subseteq [0,1]$ of length $>0$. Then
     \begin{equation}\label{eqTI} 
        g^\I(u,v)=g(v,u),\qquad\mathrm{for\ any\ } (u,v)\in [0,1]\times [0,1].
     \end{equation}
\end{theorem}

Note that, for any $Q\ge 2$, pictures like those from Figures~A and~B are
symmetric, point by point, with respect to the first diagonal, and this also
holds in the case of $\FQoI$ and $\FQeI$, for any interval $\I$ symmetric with
respect to $1/2$. If $\I$ is not centered at $1/2$, then
the corresponding pictures are never symmetric point by point,
but Theorem~\ref{TheoremI}
guarantees that their limits, as $Q\to\infty$, are symmetric. 

\smallskip

By comparing Theorem~\ref{Theorem1} above with 
Theorem~\ref{Theorem1} of~\cite{Integers},
we find that the limiting local density function $g(u,v)$ is the same
for both subsequences $\FQo$ and $\FQe$.
Therefore, Corollary 1 and Proposition 1 from~\cite{Integers}
also hold for $\FQe$.

Let $\D(0,2)$ denote the set of limit points of sequences of pairs
$\left((q'_n/Q_n,q''_n/Q_n)\right)_{n\in{\bf N}}$, 
with $Q_n\rightarrow\infty$ and
$q'_n,q''_n$ denominators of consecutive fractions in $\FQe$.

\begin{corollary}\label{Corollary1}
The set $\D(0,2)$
coincides with the quadrilateral
bounded by the lines: $y = 1,$ $x = 1,$ $2x + y = 1,$ and $2y + x = 1$.
\end{corollary}

We remark that the more uniform distribution of $\D_{500,odd}$ compared
to $\D_{500,even}$ in Figures A and B is due to a combination 
of two factors that occur
for $Q$ small. The first one is the preponderance in $\FQo$ of pairs of
type $\TT(0)$, that is, of neighbor odd fractions in $\FQ$, and the
second is the fact that the cardinality of $\FQo$ is about twice
as large as the cardinality of $\FQe$.

The next corollary provides the probability that a
neighbor pair of denominators in $\FQe$ is \emph{small}, in the sense
that their sum does not exceed $Q$.

\begin{corollary}\label{Corollary2}
The  probability that the sum of neighbor denominators of fractions
from $\FQe$ is $\le Q$ approaches $1/6$, as $Q\to\infty$.
\end{corollary}

\noindent
\begin{figure}[t]
    \begin{minipage}[t]{0.48\linewidth}
      \centering
        \includegraphics*[width=0.98\linewidth]{desen.16}
        \caption{\small The baby puzzle $\W_i$.}\label{Figure6}
     \end{minipage}
     \hfill
    \begin{minipage}[t]{0.48\textwidth}
      \centering
        \includegraphics*[width=0.98\linewidth]{desen.17}
        \caption{\small The big puzzle $\W$.}\label{Figure7}
    \end{minipage}
\end{figure}
\indent
In order to prove Theorem~\ref{Theorem1}, we begin by
presenting in Sections~\ref{SectionG} and~\ref{SectionTk} 
some geometric prerequisites, in particular the tessellation
of the Farey triangle, and then continue 
in Sections~\ref{SectionD} and~\ref{secProofT} with the study of different
types of pairs of neighbor denominators of fractions in $\FQe$. 
These conclude with
Theorem~\ref{Theorem3}, in which $g_r(u,v)$,
the local density at level $r$, is written as a sum
of a series of special quantities  assigned to the pieces of the 
tessellations of the Farey triangle.
This result plays a key role in the proof of Theorem~\ref{Theorem1} from 
Section~\ref{SectionProofT1}, where we put together all the pieces.
Remarkably, these pieces fit into a large mosaic composed by a series of
superimposed puzzles  at odd
levels (see Figures~\ref{Figure6}-\ref{Figure7}) 
recovering in the even case the same density as that obtained
in~\cite{Integers} in the odd case.
In the last section we employ a technique which makes use of estimates 
for Kloosterman sums \cite{Esterman}, \cite{Weil}, in order
to prove Theorem~\ref{TheoremI}.

\section{Some Geometry of the Farey Fractions}\label{SectionG}

We state here the basic properties of the Farey series needed for the
rest of the paper. The first one, already mentioned in the Introduction,
says that if $a'/q'<a''/q''$ are neighbor fractions in $\FQ$,
then 
    \begin{equation}\label{eq1}
        a''q'-a'q''=1\,.
    \end{equation}
Next, suppose that $a'/q'<a''/q''<a'''/q'''$ are consecutive elements
of $\FQ$. Then, the middle fraction, called the \emph{mediant}, is given by
    \begin{equation}\label{eq2}
        \frac{a''}{q''}=\frac{a'+a'''}{q'+q'''}\,.
    \end{equation}
This shows that the mediant fraction is reduced by an integer $k$, called
the {\em index} of the Farey fraction $a'/q'$, that satisfies:
    \begin{equation}\label{eq3}
k=\frac{a'+a'''}{a''}=\frac{q'+q'''}{q''}=a''q'-a'q'''
=\bigg[\frac{Q+q'}{q''}\bigg]\,.
    \end{equation}

We state a third important property in the following lemma.

\vskip 15pt
\begin{lemma}\label{LemmaA}
The positive integers $q',q''$ are denominators of neighbor fractions in $\FQ$
if and only if 
$(q',q'')\in\T_{Q}$
and $\gcd(q',q'')=1$. Also,
the pair $(q',q'')$ appears exactly once as a pair of denominators 
of consecutive
Farey fractions.
\end{lemma}

\noindent
For the proof of relations \eqref{eq1} and\eqref{eq2}, observed 
for the first time in
particular cases by Haros and Farey, we refer to 
Hardy and Wright~\cite{HW}, while for 
Lemma~\ref{LemmaA}, and further developments, see 
Hall~\cite{Hall70},~\cite{Hall94}, Hall and Tenenbaum~\cite{HT}, and
Hall and Shiu~\cite{HS}.

Now let $(q',q'',q''',\dots ,q^{(h)})$ denote a generic $h$-tuple 
of denominators of neighbor fractions in $\FQ$. Then, we see that
relation \eqref{eq3} can be employed together with Lemma~\ref{LemmaA}
to obtain a characterization of any such $h$-tuple in terms of its first
two components. Moreover, while any pair $(q',q'')$ with coprime components
$\le Q$ does appear exactly once as a pair of neighbor denominators of 
Farey fractions, the components of longer tuples must satisfy
additional conditions in order to appear together as
neighbor denominators of fractions in $\FQ$. We write these conditions
using the index.

For any positive integer $k$, we consider the convex polygon defined by
$$
\T_{_{Q,k}}:=\big\{(x,y)\colon\ 0<x,y\le Q,\ x+y>Q,\ 
                  ky\le Q+x<(k+1)y\big\}\,. 
$$
These are quadrilaterals, except for $k=1$, when $\T_{_{Q,1}}$ is a triangle.
The vertices of $\T_{_{Q,1}}$ are
$(0,Q),\; \left(\frac Q3,\frac {2Q}3\right),\; (Q,Q)$, and for any
$k\ge 2$, the vertices of $\T_{_{Q,k}}$ are
$\left(Q,\frac{2Q}k\right);$
$\left(\frac{Q(k-1)}{k+1},\frac{2Q}{k+1}\right);$
$\left(\frac{Qk}{k+2},\frac{2Q}{k+2}\right);$
$\left(Q,\frac{2Q}{k+1}\right)$.
Scaling by a factor of $Q$, for any $k\ge 1$ we get
$\T_k=\T_{_{Q,k}}/Q$, a bounded polygon inside the unit square, 
which is independent of $Q$.

On each $\T_k$ the index function, defined by
$(x,y)\mapsto \big[\frac{1+x}{y}\big]$, is locally constant. Therefore
the polygons $\T_k$ can be used to describe the
triplets of neighbor denominators of Farey fractions. We state this 
as a lemma.
\begin{lemma}\label{LemmaB}
The positive integers $q',q'',q'''$ are denominators of consecutive
fractions in $\FQ$ and $k=\frac{q'+q'''}{q''}$
if and only if $(q',q'')\in\T_{_{Q,k}}$ and $\gcd(q',q'')=1$. 
\end{lemma}

We remark that the sets $\T_{k}$, with $k\ge 1$, are disjoint and 
they form a partition of $\T$.

In many instances, one can estimate
the number of Farey fractions with a certain property by
counting the number of lattice points in a suitable domain. In this
respect, the following lemma, 
which is a variation of Lemma 2 from~[1], is very useful.
For any domain $\Omega\subset\RR^2$ we
denote:
  \begin{equation*}
    \begin{split}
          \Noo(\Omega)&:=\#\left\{(x,y)\in \Omega\cap
              \ZZ^2\colon~x~\text{odd},~y~\text{odd},\ \gcd(x,y) = 1 \right\}\,,\\
          \Neo(\Omega)&:=\#\left\{(x,y)\in \Omega\cap
              \ZZ^2\colon~x~\text{even},~y~\text{odd},\ \gcd(x,y) = 1 \right\}\,,\\
          \Noe(\Omega)&:=\#\left\{(x,y)\in \Omega\cap
              \ZZ^2\colon~x~\text{odd},~y~\text{even},\ \gcd(x,y) = 1 \right\}\,.             
    \end{split}
  \end{equation*}

These numbers can be estimated by using M\"obius
summation. One has the following asymptotic formulas.

\begin{lemma}[\cite{BCZ1}, Corollary 3.2]\label{Lemma1}
\noindent
 Let $R_1$, $R_2$ $> 0 $, and $R\ge \min(R_1,R_2)$. Then, for 
any region $\Omega\subseteq [0,R_1]\times [0,R_2]$ with rectifiable 
boundary, we have:
 \begin{eqnarray*}
  \Noo(\Omega) &=& 2\Area(\Omega)/\pi^2 + O(C_{R,\Omega})\,, \\ 
  \Noe(\Omega) &=& 2\Area(\Omega)/\pi^2 + O(C_{R,\Omega})\,,  
\end{eqnarray*}
 where $C_{R,\Omega} = \Area(\Omega)/R + R + 
\length(\partial\Omega)\log R$.
\end{lemma}
Then, for $\Omega$ as in Lemma~\ref{Lemma1}, one immediately gets 
$$  N_{\text{even},\text{odd}}(\Omega) = 2\Area(\Omega)/\pi^2 + O(C_{R,\Omega}).$$
\section{The polygons $\Tk$}\label{SectionTk}
We consider the map $T\colon\ \T\rightarrow\T$, defined by
     \begin{equation*}
       T(x,y)=\Big( y,\Big[ \tfrac{1+x}{y} \Big] y-x\Big),
     \end{equation*}
which was introduced and studied in~\cite{BCZ2}.
This transformation is invertible and its inverse is given by
     \begin{equation*}
       T^{-1}(x,y)=\Big(\Big[ \tfrac{1+y}{x} \Big]x-y,x\Big).
     \end{equation*}
One should notice that if $a'/q'<a''/q''<a'''/q'''$ are consecutive
elements in $\F_Q$, then $T(q'/Q,q''/Q)=(q''/Q,q'''/Q)$.
Then, for any ${\bf k}=(k_1,\ldots ,k_r)\in(\NN^*)^r$, we put
    \begin{equation*}
       \T_{\bf k}=\T_{k_1} \cap T^{-1} \T_{k_2} \cap \cdots
         \cap T^{-r+1} \T_{k_r} .
    \end{equation*}
We use the notational convention of dropping extra
parentheses, so for example when $r=1$ and $k\in\NN^*$, we have
$\T_{(k)}=\T_k=\big\{(x,y)\in\T\colon\
\big[\frac{1+x}{y}\big]=k\big\}$. Also, to avoid double subscripts, at
small levels $r$ we often write $k,l,\dots$, rather than $k_1,k_2,\dots$

We remark that at any level $r$, the convex polygons $\T_{\bf k}$ are
pieces of a partition of $\T$. The structure of these partitions is
studied in~\cite{Tessellation}. Here we only need the
polygons assigned to tuples $\kk$, whose components have special parities. 
These are those tuples $\kk$ with all components even, except for the first and
the last when $r\ge 2$. We call these tuples {\em admissible}, and for
each level $r$, we denote by  $\A(r)$ the set of admissible $r$-tuples.

The polygons needed in the sequel are listed in Table~\ref{Table1}.
Notice that here $\T_{\kk}$ is a quadrilateral, except for
$\kk=(1,2,4,1)$ and $\kk=(1,4,2,1)$, when it is a triangle.

\setlength{\doublerulesep}{1pt}
\scriptsize
\smallskip
\setlongtables
\begin{longtable}{|C|C|C|L|}
\caption{The polygons $\T_{\kk}$, for all admissible $\kk$.}\label{Table1}\\ \hline
\mathrm{Level} & {\normalfont \kk} & \mathrm{No.\,  of\, vertices} &
\mathrm{Vertices\, of\, } \T_{\kk}   \\ 
\hhline{|====|}
\endfirsthead
\multicolumn{4}{l}{\small\sl continued from previous page}\\ \hline
\mathrm{Level} & {\normalfont \kk} & \mathrm{No.\, of\, vertices} &
\mathrm{Vertices\, of\, } \T_{\kk}  \\ 
\hhline{|====|}
\endhead
\hline
\multicolumn{4}{r}{\small\sl continued on next page} \\ 
\endfoot
\hline
\endlastfoot
1 & k\ge 2 \text{, even } & 4 &  (\frac{k-1}{k+1}, \frac{2}{k+1});\; 
              (\frac{k}{k+2}, \frac{2}{k+2});\;
                (1,\frac{2}{k+1});\; (1, \frac 2k) \\ 
\hhline{|====|}
2 & 1,3 & 4 & (\frac 15, \frac 45);\; (\frac 27, \frac 57);\; (\frac 12, 1);\; (\frac 13, 1) \\ \hline
2 & (1,l) \text{, with $l\ge 5$ odd} & 4 & (\frac{l-3}{l+1}, \frac{l-1}{l+1});\; (\frac{l-2}{l+2},   \frac{l}{l+2});\;
        (\frac{l-1}{l+1}, 1);\; (\frac{l-2}{l}, 1) \\ \hline
2 & 3,1 & 4 & (\frac12, \frac12);\; (\frac47, \frac37);\; (1, \frac35);\; (1, \frac23)  \\ \hline
2 & (k,1) \text{, with $k\ge 5$ odd} & 4 & (\frac{k-1}{k+1}, \frac2{k+1});\; (\frac k{k+2},
\frac 2{k+2});\; (1,\frac 2{k+1});\; (1, \frac 2k)  \\ 
\hhline{|====|}
3 & 1,2,3 & 4 & (\frac 17, \frac 67);\; (\frac 15, \frac 45);\; (\frac 27, 1);\; (\frac 15, 1) \\ \hline
3 & 3,2,1 & 4 & (\frac 47, \frac 37);\; (\frac 35, \frac 25);\; (1, \frac 47);\; (1, \frac 35) \\ \hline
3 & 1,4,1 & 4 & (\frac 27, \frac 57);\; (\frac 13, \frac 23);\; (\frac 47, 1);\; (\frac 12, 1) \\ \hline
3 & (1,l,1) \text{, with $l\ge 6$, even} & 4 & (\frac{l-3}{l+1}, \frac{l-1}{l+1});\; (\frac{l-2}{l+2},
\frac{l}{l+2});\; (\frac{l-1}{l+1}, 1);\; (\frac{l-2}{l}, 1) \\
\hhline{|====|}
4 & 1,2,2,3 & 4 & (\frac 19, \frac 89);\; (\frac 17, \frac 67);\; (\frac
15, 1);\; (\frac 17, 1) \\ \hline
4 & 3,2,2,1 & 4 & (\frac 35, \frac 25);\; (\frac 57, \frac 37);\; (1,
\frac 59);\; (1, \frac 47) \\ \hline
4 & 1,2,4,1 & 3 & (\frac 15, \frac 45);\; (\frac 13, 1);\; (\frac 27, 1) \\ \hline
4 & 1,4,2,1 & 3 & (\frac 13, \frac 23);\; (\frac 35, 1);\; (\frac 47,
1) \\ 
\hhline{|====|}
r & 1,2,\dots,2,3 & 4 & (\frac{1}{2r+1}, \frac{2r}{2r+1});\; (\frac{1}{2r-1}, \frac{2r-2}{2r-1});\; (\frac{1}{2r-3}, 1);\; (\frac{1}{2r-1},1) \\ \hline
r & 3,2,\dots,2,1 & 4 & (\frac{2r-5}{2r-3}, \frac{r-2}{2r-3});\;
(\frac{2r-3}{2r-1}, \frac{r-1}{2r-1});\; (1, \frac{r+1}{2r+1});\; (1, \frac{r}{2r-1})
\end{longtable}
\normalsize

We shall denote $\T_{Q,\kk}:=Q\T_{\kk}$.

\section{The components of $\D(0,2)$}\label{SectionD}

Since any two consecutive denominators of fractions in $\FQ$ are
coprime, it follows that between any
two neighbor fractions in $\FQe$ there must be at least one
odd fraction from $\FQ$. We note that the number of
these odd fractions may be quite large as $Q$ increases. 
We classify the pairs of neighbor even fractions according
to the number of odd intermediate fractions. So, we say that the pair $(\g',\g'')$ of
consecutive fractions in $\FQe$ is of \emph{type
$\TT(r)$} if there exist exactly $r$ odd fractions between $\g'$ and
$\g''$ in $\FQ$. In this case, we shall also say that the pairs $(q',q'')$
and $(q'/Q,q''/Q)$ are of type $\TT(r)$, 
where $q',q''$ are the denominators of $\g',\g''$,
respectively.

Geometry inside the Farey triangle helps one locate the points
corresponding to such pairs.
Next, we determine, one by one, the contribution of pairs of each type
$\TT(r)$ to $\D(0,2)$.
We denote by $\E(r)$ the contribution to $\D(0,2)$ of points of order
$\TT(r)$. Then, we have
        \begin{equation}\label{eqall}
                \D(0,2)=\bigcup_{r=1}^\infty \E(r)\,.
        \end{equation}

In the following we find explicitly each set $\E(r)$.
\subsection{Points of type $\TT(1)$}
By \eqref{eq3}, it follows that pairs of fractions of this type have as
denominators the end points of a triple $(q',q'',kq''-q')$, with
$q'$ even, $q''$ odd and $k_1=k=\big[\frac{Q+q'}{q''}\big]$  even. 
This means that, for any even
$k\ge 2$, we need to retain the lattice points in the domain
  \begin{equation*}
    \begin{split}
     \U_{Q,k}=\big\{(q',kq''-q')\colon\ (q',q'')\in\T_{_{Q,k}}\big\}\,,
    \end{split}
  \end{equation*}
with $q'$ even, $q''$ odd, and $\gcd(q',q'')=1$.
In the limit, when $Q\rightarrow\infty$, these points produce a subset
of $\D(0,2)$ that is dense in the quadrilateral with vertices
   \begin{equation*}
    \begin{split}
     \U_{k}=\Big\{\Big(\frac{k-1}{k+1}, 1\Big);\  
                \Big(\frac k{k+2}, \frac k{k+2}\Big);\  
                  \Big(1, \frac{k-1}{k+1}\Big);\
                         (1, 1) 
                \Big\}\,.
    \end{split}
  \end{equation*}
Thus, we have
        \begin{equation}\label{eqlevel1}
                \E(1)=\bigcup_{\substack{k=2 \\ k\ 
        \text{even}}}^\infty\U_{k}\,.\footnote{Most often, we
denote a polygon by its vertices and, depending on the context, the
same notation will be used to indicate its topological closure.}
        \end{equation}

We remark that all $\U_k$, not only those with $k$ even, contribute to
$\D(1,2)$ (see~\cite{Integers}). 

\subsection{Points of type $\TT(2)$}

These points come from $4$-tuples of parity $(e,o,o,e)$ and this
requires both $k$ and $l$ to be odd. Let 
  \begin{equation*}
    \begin{split}
     \U_{Q,k,l}=\big\{\big(q', l(kq''-q')-q''\big)\colon\ 
                        (q',q'')\in\T_{_{Q,k,l}}\big\}\,,
    \end{split}
  \end{equation*}
and we need to pick up the lattice points in $\U_{Q,k,l}$ where
$q'$ even, $q''$ odd, and $\gcd(q',q'')=1$. 
In the limit, when
$Q\rightarrow\infty$, these points produce sequences of subsets
of $\D(0,2)$ that are dense in each of the quadrilaterals with vertices
   \begin{equation*}
    \begin{split}
     \U_{1,3}=&\big\{(1/5, 1);\; (2/7, 4/7);\; (1/2, 1/2);\; (1/3, 1)\big\}\,,\\
     \U_{3,1}=&\big\{(1/2, 1/2);\; (4/7, 2/7);\; (1, 1/5);\; (1, 1/3)\big\}\,,\\
     \U_{1,l}=&\Big\{\Big(\frac{l-3}{l+1}, 1\Big);\; 
                \Big(\frac{l-2}{l+2}, \frac{l}{l+2}\Big);\; 
                \Big(\frac{l-1}{l+1}, \frac{l-1}{l+1}\Big);\; 
                \Big(\frac{l-2}{l}, 1\Big)
        \Big\}\,,\quad \text{for $l\ge 5$}\,,\\
     \U_{k,1}=&\Big\{\Big(\frac{k-1}{k+1}, \frac{k-1}{k+1}\Big);\; 
                \Big(\frac{k}{k+2}, \frac{k-2}{k+2}\Big);\; 
                \Big(1, \frac{k-3}{k+1}\Big);\; 
                \Big(1, \frac{k-2}{k}\Big)
        \Big\}\,,\quad \text{for $k\ge 5$}\,.\\
    \end{split}
  \end{equation*}
Then
        \begin{equation}\label{eqlevel2}
                \E(2)=\U_{1,3}\cup \U_{3,1}\cup
                \bigcup_{\substack{k=5 \\ k\ \text{odd}}}^\infty\U_{k,1}
                \cup    
                \bigcup_{\substack{l=5 \\ l\ \text{odd}}}^\infty\U_{1,l}        
\,.
        \end{equation}
\subsection{Points of type $\TT(3)$}

These points are produced by $5$-tuples of parity $(e,o,o,o,e)$ and this
requires both $k$ and $m$ to be odd and $l$ even. Let 
  \begin{equation*}
    \begin{split}
     \U_{Q,k,l,m}=\big\{\big(q', m\big(l(kq''-q')-q''\big)-(kq''-q')\big)\colon\ 
                        (q',q'')\in\T_{_{Q,k,l,m}}\big\}\,,
    \end{split}
  \end{equation*}
and we need to pick up the lattice points in $\U_{Q,k,l,m}$ where
$q'$ even, $q''$ odd, and $\gcd(q',q'')=1$. 
In the limit, when
$Q\rightarrow\infty$, these points produce the sequence of subsets
of $\D(0,2)$ that are dense in each of the quadrilaterals with vertices
   \begin{equation*}
    \begin{split}
     \U_{1,2,3}=&\big\{(1/7, 1);\; (1/5, 3/5);\; (2/7, 4/7);\; (1/5, 1)\big\}\,,\\
     \U_{3,2,1}=&\big\{(4/7, 2/7);\; (3/5, 1/5);\; (1, 1/7);\; (1, 1/5)\big\}\,,\\
     \U_{1,4,1}=&\big\{(2/7, 4/7);\; (1/3, 1/3);\; (4/7, 2/7);\; (1/2, 1/2)\big\}\,,\\
     \U_{1,l,1}=&\Big\{\Big(\frac{l-3}{l+1}, \frac{l-1}{l+1}\Big);\; 
                \Big(\frac{l-2}{l+2}, \frac{l-2}{l+2}\Big);\; 
                \Big(\frac{l-1}{l+1}, \frac{l-3}{l+1}\Big);\; 
                \Big(\frac{l-2}{l}, \frac{l-2}{l}\Big)
        \Big\}\,,\quad \text{for $l\ge 6$}\,.
    \end{split}
  \end{equation*}
Then
        \begin{equation}\label{eqlevel3}
                \E(3)=\U_{1,2,3}\cup \U_{3,2,1}\cup 
                \bigcup_{\substack{l=4 \\ l\ \text{even}}}^\infty\U_{1,l,1}\,.
        \end{equation}
\subsection{Points of type $\TT(4)$}

These points are produced by $6$-tuples of parity $(e,o,o,o,o,e)$ and this
requires both $k,n$ to be odd and $m,l$ even. Let 
\begin{multline*}
     \U_{Q,k,l,m,n}=\big\{\big(q',
                n\big(m\big(l(kq''-q')-q''\big)-(kq''-q')\big)-
                                \big(l(kq''-q')-q''\big)\colon\ \\ 
                        (q',q'')\in\T_{_{Q,k,l,m,n}}\big\}\,.
\end{multline*}
As before, we need to pick up the lattice points
in $\U_{Q,k,l,m,n}$  with $q'$ even, $q''$ odd, and $\gcd(q',q'')=1$. 
By Table~\ref{Table1}, we know that there are only four such sets. 
In the limit, as $Q\rightarrow\infty$, we obtain two 
quadrilaterals and two triangles: 
   \begin{equation*}
    \begin{split}
     \U_{1,2,2,3}=&\big\{(1/9, 1);\; (1/7, 5/7);\; (1/5, 3/5);\; (1/7, 1)\big\}\,,\\
     \U_{3,2,2,1}=&\big\{(3/5, 1/5);\; (5/7, 1/7);\; (1, 1/9);\; (1, 1/7)\big\}\,,\\
     \U_{1,2,4,1}=&\big\{(1/5, 3/5);\; (1/3, 1/3);\; (2/7, 4/7)\big\}\,,\\
     \U_{1,4,2,1}=&\big\{(1/3, 1/3);\; (3/5, 1/5);\; (4/7, 2/7)\big\}\,.\\
    \end{split}
  \end{equation*}
Then
        \begin{equation}\label{eqlevel4}
                \E(4)=\U_{1,2,2,3}\cup \U_{3,2,2,1}\cup 
                      \U_{1,2,4,1}\cup \U_{1,4,2,1}\,.
        \end{equation}

\subsection{Points of type $\TT(r),\ r\ge 5$}

These points are produced by the $(r+2)$-tuples $(q',\dots,q^{(r+2)})$ 
of parity $(e,o,\dots,o,e)$.
It turns out that the only possible corresponding $r$-tuples $\kk$ are
$(1,2,\dots,2,3)$ and its symmetric $(3,2,\dots,2,1)$, with $(r-2)$
$2$'s between the endpoints.
This allows us to find a closed formula for
$q^{(r+2)}$. We find that
$q^{(r+2)}=-(2 r-1)q' + 2 q''$ in the first case,
and $q^{(r+2)}=-q' + 2 q''$ in the second case. Thus, we define
  \begin{equation*}
    \begin{split}
     \U_{Q,1,2,\dots,2,3}=\big\{\big(q', -(2 r-1)q' + 2 q''\big)\colon\ 
                        (q',q'')\in\T_{_{Q,1,2,\dots,2,3}}\big\}\,,
    \end{split}
  \end{equation*}
and
  \begin{equation*}
    \begin{split}
     \U_{Q,3,2,\dots,2,1}=\big\{\big(q', -q' + 2 q''\big)\colon\ 
                        (q',q'')\in\T_{_{Q,3,2,\dots,2,1}}\big\}\,.
    \end{split}
  \end{equation*}

As before, only the lattice points with $q'$ even, $q''$ odd, and
$\gcd(q',q'')=1$ should be considered, and in the limit, when 
$Q\rightarrow\infty$, for any given $r$, we obtain the quadrilaterals:
   \begin{equation*}
    \begin{split}
             \U_{1,2,\dots,2,3}=&\Big\{\Big(\frac{1}{2r+1}, 1\Big);\; 
                \Big(\frac{1}{2r-1}, \frac{2r-3}{2r-1}\Big);\; 
                \Big(\frac{1}{2r-3}, \frac{2r-5}{2r-3}\Big);\; 
                \Big(\frac{1}{2r-1}, 1\Big)
        \Big\}\,,\quad \text{for $r\ge 5$}\,,\\
             \U_{3,2,\dots,2,1}=&\Big\{\Big(\frac{2r-5}{2r-3}, \frac{1}{2r-3}\Big);\; 
                \Big(\frac{2r-3}{2r-1}, \frac{1}{2r-1}\Big);\; 
                \Big(1, \frac{1}{2r+1}\Big);\; 
                \Big(1,\frac{1}{2r-1}\Big)
        \Big\}\,,\quad \text{for $r\ge 5$}\,.
    \end{split}
  \end{equation*}
Then
        \begin{equation}\label{eqlevelr}
                \E(r)=\U_{1,2,\dots,2,3}\cup \U_{3,2,\dots,2,1}\,,
                \quad\text{ for $r\ge 5$}.
        \end{equation}


In Figure~\ref{Figure4}, 
one can see a representation of $\D(0,2)$ covered by
$\U_*$. In addition, the union from the
right hand side of \eqref{eqall} gives a first
hint on the local densities on $\D(0,2)$. Complete calculations
are postponed to Section~\ref{secProofT}.

\section{The density of points of type $\TT(r)$}\label{secProofT}
Let $g(x,y)$ be the function that gives the local density of 
points $(q'/Q,q''/Q)$ in the unit square as $Q\to\infty$, 
and let $g_r(x,y)$ be the local density 
in the unit square of the points $(q'/Q,q''/Q)$ of type $\TT(r)$, as
$Q\rightarrow\infty$. At any point $(u,v)\in(0,1)^2$, 
this local density is defined by 
    \begin{equation}\label{eqdefgr}
      g_r(u,v):=\lim_{\Area(\Box)\to 0}
           \frac{\lim\limits_{Q\to \infty}
             \frac{\#\big(\Box\cap\D_Q(0,2)/Q\big)}{\#\D_Q(0,2)}}
                  {\Area(\Box)}\,,
    \end{equation}
where $\Box\subset\RR^2$ are squares centered at $(u,v)$.
They are related by 
  \begin{equation}\label{eqgg}
    g(u,v)=\sum_{r=1}^\infty g_r(u,v)\,,
  \end{equation}
provided we show that each local density $g_r(u,v)$ exists as $Q\rightarrow\infty$.
In the following we find each $g_r(u,v)$. 
The proof generalizes that of Theorem 1 from~\cite[Section 3.2]{Integers}.

\subsection{Generalities on $g_r$}
Let  $r\ge 1$, let $(x_0,y_0)$ be a fixed point inside the unit square $(0,1)^2$, and 
fix a small $\eta>0$.
In the following, the superscript $\LL$ indicates the last element of a
tuple. For example, $q^{\LL}=q^\LL_r(q',q'')$ is the 
denominator of the $(r+2)$-nd fraction in $\FQ$, starting with $q', q''$.
We denote the square centered at $(x_0,y_0)$ by 
$\Box=\Box_\eta(x_0,y_0)=(x_0-\eta,x_0+\eta)\times(y_0-\eta,y_0+\eta)$.

By definition,
any pair $(q',q^\LL)$ of type $\TT(r)$ is generated by an $(r+2)$-tuple
$(q',q'',\dots,q^\LL)$
of denominators of consecutive fractions in $\FQ$, 
with $q',q^\LL$ even, and the rest of
the components odd. We consider the set $\B_Q$ of pairs $(q',q'')$
of type $\TT(r)$ for which the corresponding point
$(q'/Q,q^\LL/Q)$ falls in $\Box$, that is,
  \begin{equation*}
     \B_Q(r)=\left\{(q',q'')\in \NN^2\colon\ 
        \begin{array}{l}  1\le q',q''\le Q,\ \gcd(q',q'') = 1,\ q'+q''>Q,
                 \\ \displaystyle
        q'~ \text{even},\ q''~ \text{odd}; \ \kk(q',q'')\in\A(r),\ 
              (q',q^\LL(r))\in Q\cdot \Box
        \end{array} \right\}.
  \end{equation*}  
The cardinality of $\B_Q$ is $\#\B_Q=N_{\text{even},\text{odd}}(\Omega_Q) $, where
$\Omega_Q=\Omega_Q(r)=\Omega_Q(x_0,y_0,\eta)(r)$ is given by 
  \begin{equation*}
     \Omega_Q(r)=\left\{(x,y)\in \RR^2\colon\ 
        \begin{array}{l}  
              1\le x,y\le Q,\  x+y>Q,\  \\ \displaystyle
         \kk(x,y)\in\A(r),\     (x,x_r^\LL(x,y))\in Q\cdot\Box 
        \end{array} 
        \right\}.
  \end{equation*}
Here $x^\LL=x^\LL_r(x,y)$ is the $(r+2)$-nd element of the sequence defined
recursively by: $x_{-1}=x,\ x_0=y$ and $x_j=k_jx_{j-1}-x_{j-2}$, where
$k_j=\big[\frac{1+x_{j-2}}{x_{j-1}}\big]$, for $1\le j\le r$ and 
$\kk(x,y)=(k_1,\dots,k_r)\in\A(r)$.

Multiplying by $1/Q$, we obtain the bounded set
  \begin{equation*}
     \Omega(r)=\left\{(x,y)\in (0,1)^2\colon\ 
                x+y>1,\ \kk(x,y)\in\A(r),\ 
                (x,x^\LL_r(x,y))\in \Box\right\},
  \end{equation*}
and $Q\cdot \Omega(r)=\Omega_Q(r)$. We are interested in the
area of $\Omega(r)$,
since, by Lemma~\ref{Lemma1}, we know that
    \begin{equation}\label{eqBQ}
       \#\B_Q(r)=\frac{2Q^2 \Area\big(\Omega(r)\big)}{\pi^2}
                      +O(Q\log Q)\,.
    \end{equation}

Next we split $\Omega(r)$ into the pieces given by the shadows left on
it by each of the sets
    \begin{equation*}
       \U_\kk(r):=\left\{(x,y)\in (0,1)^2\colon\ \kk_r(x,y)=\kk \right\},
        \quad \text{for\ $\kk\in\A(r)$}\,.
    \end{equation*}
Since $\Omega(r)\subset\T$ and $\U_\kk(r)\cap \T=\T_\kk$, it follows that
$\Omega(r)\cap\U_\kk(r)=\T_\kk\cap P_\kk(r)$, where
$\P_{\kk}(r)=\P_{\kk,r}(x_0,y_0,\eta)$ is the parallelogram
   \begin{equation*}
            \P_\kk(r)=\left\{(x,y)\in \RR^2\colon\ 
                        (x,x^\LL_r(x,y))\in\Box_\eta(x_0,y_0)
                      \right\}.
    \end{equation*}
Thus, we have obtained 
     \begin{equation}\label{eqAO}
         \Area(\Omega(r))=\sum_{\kk\in\A(r)} \Area(\T_\kk\cap \P_\kk(r))\,.
     \end{equation}
A compactness argument shows that, although $\A(r)$ may be infinite,
only finitely many terms of the series
are non-zero. Our next objective is to make explicit
their size in terms of the position of $(x_0,y_0)$.
But first we need to establish a concrete expression for $\P_\kk(r)$.

\subsection{The index $p_r(\kk)$ and the parallelogram $\P_\kk(r)$}
First, we define the sequence of polynomials $p_r(\kk)$ by: 
$p_0(\cdot)=1$, $p_1(k_1)=k_1$, and then recursively, for any
$r\ge 2$,
    \begin{equation}\label{eqR}
   p_r(k_1,\dots,k_r)=k_rp_{r-1}(k_1,\dots,k_{r-1})-p_{r-2}(k_1,\dots,k_{r-2}).
    \end{equation}
The first polynomials are:
        \begin{align*}
        p_2(\kk)&=k_1k_2-1;\\
        p_3(\kk)&=k_1k_2k_3-k_1-k_3;\\
        p_4(\kk)&=k_1k_2k_3k_4-k_1k_2-k_1k_4-k_3k_4+1;\\
        p_5(\kk)&=k_1k_2k_3k_4k_5-k_1k_2k_3-k_1k_2k_5-k_1k_4k_5-k_3k_4k_5+k_1+k_3+k_5.
        \end{align*}  
Some other particular values we need later are:
  \begin{equation}\label{eqPartVal}
    \begin{split}
        p_r(1,2,\dots,2)\phantom{,3} &=1,\phantom{r+1}\quad\text{for $r\ge 2$,}\\
        p_r(1,2,\dots,2,3)&=2,\phantom{r+1}\quad\text{for $r\ge 2$,}\\
        p_r(2,\dots,2,3)&=2r+1,\quad\text{for $r\ge 2$,}
    \end{split}
  \end{equation}
Also, we remark the symmetry property: 
        \begin{equation}\label{eqSym}
                p_r(k_r,\dots,k_1)=p_r(k_1,\dots,k_r)\,.
        \end{equation}
More on this fundamental sequence of polynomials 
can be found in~\cite{Tessellation}.

Next, let us notice that $x_r^\LL(x,y)$ is a linear combination of $x$
and $y$:
  \begin{equation}\label{eqLinComb}
    \begin{split}
        x_r^\LL(x,y)=p_r(k_1,\dots,k_r)y-p_{r-1}(k_2,\dots,k_r)x\,.
    \end{split}
  \end{equation}

Returning now to our parallelogram, by \eqref{eqLinComb} 
we see that this is the set of
points $(x,y)\in\RR^2$ that satisfy the conditions:
  \begin{equation*}
    \begin{split}
        \begin{cases}
                x_0-\eta<x<x_0+\eta,\\
        y_0-\eta<p_r(k_1,\dots,k_r)y-p_{r-1}(k_2,\dots,k_r)x<y_0+\eta\,.
        \end{cases}
    \end{split}
  \end{equation*}

From this, we see that the area of $\P_\kk(r)$ is
  \begin{equation}\label{eqAreaP}
    \begin{split}
        \Area(\P_\kk(r))=\frac{4\eta^2}{p_r(\kk)}\,,
    \end{split}
  \end{equation}
and its center has coordinates
  \begin{equation}\label{eqCenter}
    \begin{split}
        C_\kk(r)=\left(x_0,\;\frac{p_{r-1}(k_2,\dots,k_r)}{p_r(k_1,\dots,k_r)}x_0
                +\frac{1}{p_r(k_1,\dots,k_r)}y_0\right).
    \end{split}
  \end{equation}

\subsection{The density $g_r(u,v)$ II}

In order to obtain a concrete expression for the density, one requires
an explicit form of the series in \eqref{eqAO}. 
Since $\eta>0$ can be chosen as small as we please, it follows that
the summands there depend on the position of $(x_0,y_0)$ with respect
to $\T_\kk$. 
In the following, we shall assume that $\eta$ is small enough. We may
also assume that $\kk$ is bounded, since all the parallelograms $\P_\kk(r)$
are contained in the vertical strip given by the
inequalities $x_0-\eta<x<x_0+\eta$, and since only finitely many
polygons $\T_\kk$ intersect this strip. 

Let now $\kk$ be an admissible $r$-tuple. We check what happens when
$C_\kk(r)$ lies on the interior $\overset{\circ}{\T_\kk}$
of $\T_\kk$, on the edges of $\T_\kk$,
or in the set $V(\kk)$ of vertices of $\T_\kk$. 

Firstly, if $C_\kk(r)\in\overset{\circ}{\T_\kk}$, 
it follows that $\P_\kk(r)\subset\T_\kk$, 
so $\Area(\T_\kk\cap \P_\kk(r))=\Area(\P_\kk(r))$. 
Secondly, if $C_\kk(r)\in\partial{\T_\kk}\setminus V(\T_\kk)$, then 
$\Area(\T_\kk\cap \P_\kk(r))=\Area(\P_\kk(r))/2$, since 
any line that crosses $\P_\kk(r)$  through its center cuts 
the parallelogram into two pieces of equal area.
Thirdly, suppose $C_\kk(r)\in V(\T_\kk)$. Then $\Area(\T_\kk\cap\P_\kk(r))$
depends on the angle formed by the corresponding
edges, and these angles may differ for
different vertices of $\T_\kk$ or for different values of $\kk$.  
We have collected all the results in Table~\ref{Table2}. 
In the calculations, we have made use of the relations
\eqref{eqPartVal}, \eqref{eqSym} and \eqref{eqLinComb}.
For each $V\in V(\T_\kk)$ and $C_\kk(r)=V$, for a concise
presentation, we have translated 
$\T_\kk\cap\P_\kk(r)$ with a vector $V$ to the origin. Also, we remark
that the size of $\T_\kk\cap\P_\kk(r)$ is always
proportional to $\eta$. It follows that
$\alpha_\kk(V):=\Area(\T_\kk\cap\P_\kk(r))/\eta^2$ is independent 
of $\eta$.

The nice thing about this calculation is that it has an
error-correcting check. This is due to another remarkable property of
the parallelogram,  which implies that
the entries in the column of $\alpha_\kk(V)$ satisfy the relations
  \begin{equation}\label{eqremarq}
    \begin{split}
        \sum_{V}\alpha_\kk(V)=
        \begin{cases}
                \frac{\Area(\P_\kk(r))}{\eta^2}= \frac{4}{p_r(\kk)},
                        \quad\text{if $\T_\kk$ is a quadrilateral;}\\
                \frac{\Area(\P_\kk(r))}{2\eta^2}=\frac{2}{p_r(\kk)}, 
                        \quad\text{if $\T_\kk$ is a triangle,}
        \end{cases}
    \end{split}
  \end{equation}
for each admissible $\kk$. Here the summation is over all the vertices
of $\T_{\kk}$. For example, if $\kk=(1,l,1)$, with $l\ge 6$ even, we have
  \begin{equation*}
        \frac{2l-1}{2(l-1)l}+\frac{l+2}{(l-2)l}
        +\frac{2l-1}{2(l-1)l}+\frac{l}{(l-2)(l-1)}=\frac{4}{l-2}\,
  \end{equation*}
and for $\kk=(1,2,4,1)$, we have
  \begin{equation*}
        \frac{7}{24}+\frac{3}{40}+\frac{19}{30}=1\,.
  \end{equation*}

\setlength{\doublerulesep}{1pt}
\scriptsize
\smallskip
\setlongtables
\begin{longtable}{|C|C|L|C|}
\caption{The vertices and the area of the polygons $\T_{\kk}\cap P_{\kk}(r)$
for $C_{\kk}(r)$ equal to each of the vertices of $\T_{\kk}$. Each
polygon was scaled by $1/\eta$ and translated into the origin. 
}\label{Table2}\\ \hline
{\normalfont \kk} &\text{Vertex }  V=C_{\kk} & 
\mathrm{Vertices\ of\ } \frac {1}{\eta}(\T_{\kk}\cap P_{\kk}(r))-C_{\kk}(r) &
\alpha_\kk(V)
\\ 
\hhline{|====|}
\endfirsthead
\multicolumn{4}{l}{\small\sl continued from previous page}\\ \hline
{\normalfont \kk} &\text{Vertex }  V=C_{\kk}(r) & 
\mathrm{Vertices\ of\ } \frac {1}{\eta}(\T_{\kk}\cap P_{\kk}(r))-C_{\kk}(r) &
\alpha_\kk(V)
\\ 
\hhline{|====|}
\endhead
\hline
\multicolumn{4}{r}{\small\sl continued on next page} \\ 
\endfoot
\hline
\endlastfoot
k\ge 2, \text{ even}  
& (\frac{k-1}{k+1}, \frac{2}{k+1} )
& (0, 0);\; (\frac{1}{k+1}, -\frac{1}{k+1});\; (1, 0);\; (1, \frac 1k)
& \frac{2k+1}{2k(k+1)} \\ \hline
& (\frac{k}{k+2}, \frac{2}{k+2})
& (0, 0);\; (1, \frac{1}{k+1});\; (1, \frac{2}{k});\; (-\frac{1}{k+1}, \frac{1}{k+1})
& \frac{k+2}{k(k+1)} \\ \hline
&  (1, \frac{2}{k+1})
 & (0, 0);\; (0, \frac 1k);\; (-1, 0);\; (-1, -\frac{1}{k+1})
& \frac{2k+1}{2k(k+1)}  \\ \hline
& (1, \frac 2k)
& (0, 0);\; (-1, -\frac 1k);\; (-1, -\frac 2k);\; (0, -\frac 1k) 
& \frac 1k \\ 
\hhline{|====|}
1,3
& (\frac 15, \frac 45) 
& (0, 0);\; (\frac 15, -\frac 15);\; (1,1);\; (1, \frac{3}{2}) 
& \frac{9}{20} \\ \hline
& (\frac 27, \frac 57) 
& (0, 0);\; (1, \frac 43);\; (1, 2);\; (-\frac 15, \frac 15) 
& \frac{19}{30} \\ \hline
& (\frac 12, 1) 
& (0, 0);\; (-\frac 13, 0);\; (-1, -1);\; (-1, -\frac 43)
& \frac 13 \\ \hline
& (\frac 13, 1) 
& (0, 0);\; (-1, -\frac 32);\; (-1, -2);\; (\frac 13, 0)
& \frac{7}{12} \\ 
\hhline{|====|}
(1,l), \text{ with $l\ge 5$ odd}
 & (\frac{l-3}{l+1}, \frac{l-1}{l+1})
& (0, 0);\; (\frac{2}{l+1}, \frac{1}{l+1});\; (1, 1);\; (1, \frac{l}{l-1})
& \frac{l}{(l-1)(l+1)} \\ \hline
& (\frac{l-2}{l+2}, \frac{l}{l+2})
& (0, 0);\; (1, \frac{l+1}{l});\; (1, \frac{l+1}{l-1});\; 
(-\frac{2}{l+1}, -\frac{1}{l+1})
& \frac{2l^2+5l+1}{2(l-1)l(l+1)} \\ \hline
& (\frac{l-1}{l+1}, 1)
& (0, 0);\; (-\frac{1}{l}, 0);\; (-1, -1);\; (-1, -\frac{l+1}{l})
& \frac{1}{l} \\ \hline
& (\frac{l-2}{l}, 1)
& (0, 0);\; (-1, -\frac{l}{l-1});\; (-1, -\frac{l+1}{l-1});\; (\frac{1}{l}, 0) 
& \frac{2l+1}{2l(l-1)} \\
\hhline{|====|}
3,1
& (\frac 12, \frac 12)
& (0, 0);\; (\frac 13, -\frac 13);\; (1, 0);\; (1, \frac 13)
& \frac 13  \\ \hline
& (\frac 47, \frac 37)
& (0, 0);\; (1, \frac 25);\; (1, 1);\; (-\frac 13, \frac 13)
& \frac{19}{30} \\ \hline
& (1, \frac 35 )
& (0, 0);\; (0, \frac 12);\; (-1, 0);\; (-1, -\frac 25)
& \frac{9}{20} \\ \hline
& (1, \frac 23 )
& (0, 0);\; (-1, -\frac 13);\; (-1, -1);\; (0, -\frac 12)
& \frac{7}{12} \\ 
\hhline{|====|}
(k,1), \text{ with $k\ge 5$ odd}
& (\frac{k-1}{k+1}, \frac{2}{k+1}) 
& (0, 0);\; (\frac 1k, -\frac 1k);\; (1, 0);\; (1, \frac 1k)
& \frac 1k \\ \hline
& (\frac{k}{k+2}, \frac{2}{k+2})
& (0, 0);\; (1, \frac{1}{k+1});\; (1, \frac{2}{k-1});\; (-\frac 1k, \frac 1k)
& \frac{2k^2+5k+1}{2(k-1)k(k+1)}  \\ \hline
& (1, \frac{2}{k+1})
& (0, 0);\; (0, \frac{1}{k-1});\; (-1, 0);\; (-1, -\frac{1}{k+1})
& \frac{k}{(k-1)(k+1)} \\ \hline
& (1, \frac 2k)
& (0, 0);\; (-1, -\frac 1k);\; (-1, -\frac{2}{k-1});\; (0, -\frac{1}{k-1})
& \frac{2k+1}{2k(k+1)} \\ 
\hhline{|====|}
1,2,3
& (\frac 17, \frac 67)
& (0, 0);\; (\frac 17, -\frac 17);\; (1, 2);\; (1, \frac 52)
& \frac{13}{28} \\ \hline
& (\frac 15, \frac 45)
& (0, 0);\; (1, \frac 73);\; (1, 3);\; (-\frac 17, \frac 17)
& \frac{13}{21} \\ \hline
& (\frac 27, 1) 
& (0, 0);\; (-\frac 15, 0);\; (-1, -2);\; (-1, -\frac 73)
& \frac{11}{30} \\ \hline
& (\frac 15, 1)
& (0, 0);\; (-1, -\frac 52);\; (-1, -3);\; (\frac 15, 0)
& \frac{11}{20} \\ 
\hhline{|====|}
3,2,1
& (\frac 47, \frac 37)
& (0, 0);\; (\frac 13, -\frac 13);\; (1, 0);\; (1, \frac 25) 
& \frac{11}{30}  \\ \hline
& (\frac 35, \frac 25)
& (0, 0);\; (1, \frac 37);\; (1, 1);\; (-\frac 13, \frac 13) 
& \frac{13}{21}  \\ \hline
& (1, \frac 47)
& (0, 0);\; (0, \frac 12);\; (-1, 0);\; (-1, -\frac 37) 
& \frac{13}{28}  \\ \hline
& (1, \frac 35)
& (0, 0);\; (-1, -\frac 25);\; (-1, -1);\; (0, -\frac 12) 
& \frac{11}{20} \\ 
\hhline{|====|}
1,4,1 
& (\frac 27, \frac 57)
& (0, 0);\; (\frac 15, -\frac 15);\; (1, 1);\; (1, \frac 43)
& \frac{11}{30}   \\ \hline
& (\frac 13, \frac 23)
& (0, 0);\; (1, \frac 75);\; (1, 2);\; (-\frac 15, \frac 15)
& \frac 35   \\ \hline
& (\frac 47, 1)
& (0, 0);\; (-\frac 13, 0);\; (-1, -1);\; (-1, -\frac 75)
& \frac{11}{30}  \\ \hline
& (\frac 12, 1)
& (0, 0);\; (-1, -\frac 43);\; (-1, -2);\; (\frac 13, 0)
& \frac 23 \\ 
\hhline{|====|}
(1,l,1), \text{ with $l\ge 6$ even }
& (\frac{l-3}{l+1}, \frac{l-1}{l+1})
& (0, 0);\; (\frac 2l, \frac 1l);\; (1, 1);\; (1, \frac{l}{l-1}) 
& \frac{2l-1}{2(l-1)l}  \\ \hline
&  (\frac{l-2}{l+2}, \frac{l}{l+2})
& (0, 0);\; (1, \frac{l+1}{l});\; (1, \frac{l}{l-2});\; (-\frac{2}{l}, -\frac 1l) 
& \frac{l+2}{(l-2)l}  \\ \hline
& (\frac{l-1}{l+1}, 1)
& (0, 0);\; (-\frac{1}{l-1}, 0);\; (-1, -1);\; (-1, -\frac{l+1}{l})
& \frac{2l-1}{2(l-1)l}  \\ \hline
& (\frac{l-2}{l}, 1)
& (0, 0);\; (-1, -\frac{l}{l-1});\; (-1, -\frac{l}{l-2});\; (\frac{1}{l-1}, 0) 
& \frac{l}{(l-2)(l-1)} \\ 
\hhline{|====|}
1,2,2,3
& (\frac 19, \frac 89)
& (0, 0);\; (\frac 19, -\frac 19);\; (1, 3);\; (1, \frac 72) 
& \frac{17}{36}  \\ \hline
& (\frac 17, \frac 67)
& (0, 0);\; (\frac 12, \frac 54);\; (1, 3);\; (-\frac 19, \frac 19) 
& \frac{65}{72} \\ \hline
& (\frac 15, 1) 
& (0, 0);\; (-\frac 17, 0);\; (-\frac 12, -\frac 54)
& \frac{5}{56}  \\ \hline
& (\frac 17, 1)
& (0, 0);\; (-1, -\frac 72);\; (-1, -4);\; (\frac 17, 0) 
& \frac{15}{28} \\ 
\hhline{|====|}
1,2,4,1
& (\frac 15, \frac 45) 
& (0, 0);\; (\frac 12, \frac 34);\; (1, 2);\; (1, \frac 73) 
& \frac {7}{24}  \\ \hline
& (\frac 13, 1)
& (0, 0);\; (-\frac 15, 0);\; (-\frac 12, -\frac 34) 
& \frac{3}{40}  \\ \hline
& (\frac{2}{7}, 1)
& (0, 0);\; (-1, -\frac 73);\; (-1, -3);\; (\frac 15, 0) 
& \frac{19}{30} \\ 
\hhline{|====|}
1,4,2,1
& (\frac 13, \frac 23)
& (0, 0);\; (1, \frac 54);\; (1, \frac 75) 
& \frac{3}{40} \\ \hline
& (\frac 35, 1)
& (0, 0);\; (-\frac 13, 0);\; (-1, -1);\; (-1, -\frac 54) 
& \frac{7}{24} \\ \hline
& (\frac 47, 1)
& (0, 0);\; (-1, -\frac 75);\; (-1, -2);\; (\frac 13, 0)
& \frac{19}{30} \\ 
\hhline{|====|}
3,2,2,1
& (\frac 35, \frac 25)
& (0, 0);\; (1, \frac 14);\; (1, \frac 37) 
& \frac{5}{56}  \\ \hline
& (\frac 57, \frac 37)
& (0, 0);\; (1, \frac 49);\; (1, 1);\; (-1, 0);\; (-1, -\frac 14) 
& \frac{65}{72}  \\ \hline
& (1, \frac 59)
& (0, 0);\; (0, \frac 12);\; (-1, 0);\; (-1, -\frac 49) 
& \frac{17}{36}  \\ \hline
& (1, \frac 47)
& (0, 0);\; (-1, -\frac 37);\; (-1, -1);\; (0, -\frac 12) 
& \frac{15}{28} \\ 
\hhline{|====|}
1,2,\dots,2,3
& (\frac{1}{2r+1}, \frac{2r}{2r+1})
& (0, 0);\; (\frac{1}{2r+1}, -\frac{1}{2r+1});\; (1, r-1);\; (1, \frac{2r-1}{2}) 
& \frac{4r+1}{4(2r+1)}  \\ \hline
& (\frac{1}{2r-1}, \frac{2(r-1)}{2r-1})
& (0, 0);\; (\frac 12, \frac{2r-3}{4});\; (1, r-1);\; (1, r);\; 
      (-\frac{1}{2r+1}, \frac{1}{2r+1}) 
& \frac{14r+9}{8(2r+1)}  \\ \hline
& (\frac{1}{2r-3}, 1)
& (0, 0);\; (-\frac{1}{2r-1}, 0);\; (-\frac 12, -\frac{2r-3}{4}) 
& \frac{2r-3}{8(2r-1)}  \\ \hline
& (\frac{1}{2r-1}, 1)
& (0, 0);\; (-1, -\frac{2r-1}{2});\; (-1, -r);\; (\frac{1}{2r-1}, 0) 
& \frac{4r-1}{4(2r-1)} \\ 
\hhline{|====|}
3,2,\dots,2,1
&  (\frac{2r-5}{2r-3}, \frac{r-2}{2r-3} )
& (0, 0);\; (1, \frac 14);\; (1, \frac{r-1}{2r-1}) 
& \frac{2r-3}{8(2r-1)} \\ \hline
& (\frac{2r-3}{2r-1}, \frac{r-1}{2r-1} )
& (0, 0);\; (1, \frac{r}{2r+1});\; (1, 1);\; (-1, 0);\; (-1, -\frac 14) 
& \frac{14r+9}{8(2r+1)}  \\ \hline
& (1, \frac{r+1}{2r+1})
& (0, 0);\; (0, \frac 12);\; (-1, 0);\; (-1, -\frac{r}{2r+1}) 
& \frac{4r+1}{4(2r+1)}  \\ \hline
& (1, \frac{r}{2r-1})
& (0, 0);\; (-1, -\frac{r-1}{2r-1});\; (-1, -1);\; (0, -\frac 12) 
& \frac{4r-1}{4(2r-1)} \\ \hline
\end{longtable}
\normalsize

These observations combined with \eqref{eqAreaP} put \eqref{eqAO}
into the form
     \begin{equation}\label{eqAO1}
         \Area(\Omega(r))=
                4\eta^2\sum_{C_\kk(r)\in\overset{\circ}{\T_\kk}} \frac 1{p_r(\kk)}
             +2\eta^2\sum_{C_\kk(r)\in\partial{\T_\kk}\setminus
                V(\T_\kk)} \frac 1{p_r(\kk)}
         +\eta^2\sum_{C_\kk(r)\in V(\T_\kk)} \alpha_\kk(r)\,.
     \end{equation}

We now immediately obtain a corresponding expression
for $g_r$. An application of Lemma~\ref{Lemma1} (see also
Lemma~\ref{LemmaFQeI} below)  provides
$$ 
\#\FQe=Q^2/\pi^2+O(Q\log Q). 
$$
Using \eqref{eqBQ} 
and the fact that the number of points 
$(q'/Q,q''/Q)$ from $(0,1)^2$, where $q'$ and $q''$ are 
denominators of two consecutive elements from $\FQe$ is $\#\FQe-1$, we have:
   \begin{equation}\label{eqRaportul}
        \iint\limits_{\Box_\eta(x_0,y_0)}g_r(x,y)\,dxdy
         =\lim_{Q\rightarrow \infty}
           \frac{\#\B_Q(r)}{\#\FQe-1} = 2\Area(\Omega(r))\,.
   \end{equation}
Then, by the Lesbegue differentiation theorem, we have
$ g_r(x_0,y_0) = \lim_{\eta\rightarrow 0} 
                \Area(\Omega(r))/(4\eta^2)$,
which combined with  \eqref{eqAO1} gives the following result.
\begin{theorem}\label{Theorem3}
For $(x_0,y_0)\in [0,1]^2$ and any integer $r\ge 1$, we have:
     \begin{equation}\label{eqgr} 
     g_r(x_0,y_0) =
                \sum_{C_\kk(r)\in\overset{\circ}{\T_\kk}} \frac 2{p_r(\kk)}
             +\sum_{C_\kk(r)\in\partial{\T_\kk}\setminus
                V(\T_\kk)} \frac 1{p_r(\kk)}
         +\frac 12\sum_{C_\kk(r)\in V(\T_\kk)} \alpha_\kk(r)\,,
     \end{equation}
where the sums run over $r-$tuples $\kk$ which are admissible.
\end{theorem}

We call the generic term in the first sum the {\em kernel} of
$\kk$. Notice that if $C_\kk(r)$ is on the boundary the terms added in
the second sum are equal to half of the kernel, and in the third sum the
kernel distributes in sizes proportional to the angles of different
vertices of the polygon $\T_\kk$. 

\section{Proof of Theorem~\ref{Theorem1}}\label{SectionProofT1}

Our starting point in the proof of Theorem~\ref{Theorem1} is formula
\eqref{eqgr}, in which we need to make explicit the conditions of
summation in terms of the variables. Since some of the  polygons
with all components of $\kk$ small do not follow the general pattern,
at this point we can not get an explicit closed formula for the
density, even if we restrict to pairs of a given order. Consequently,
we further split the sums on the right hand side of \eqref{eqgr}.
For this, we collect in $\huu(x_0,y_0)$ the sum of the terms from the
right hand side of \eqref{eqgr} with $r=1$ and $\kk=k\ge 4$,
even. Similarly, for $r=2$, $\hdu(x_0,y_0)$ will be the sum over
$\kk=(1,l)$ and $\kk=(k,1)$ with $k,l\ge 5$, both odd, while for
$r=3$, $\htu(x_0,y_0)$ will be the sum over $\kk=(1,l,1)$, with $l\ge
6$, even. The remaining terms will be collected separately in
$\hud(x_0,y_0)$, $\hdd(x_0,y_0)$ and $\htd(x_0,y_0)$. Thus, we have:
  \begin{equation*}
    \begin{split}
  g_j(x_0,y_0)=&h_j^\mathrm{d}(x_0,y_0)+h_j^\mathrm{u}(x_0,y_0),\quad\text{ for $j=1,2,3.$}\\
    \end{split}
  \end{equation*}

Next, we treat one by one each of these terms. We shall assume 
everywhere, unless otherwise specified, that $x_0,y_0>0$.
To shorten the notation, we shall use the variables $z_0$ and
$\z_0$ for either $x_0$ or $y_0$, with the meaning explained
in the Introduction for $z,\z$ and $u,v$. In the
same way, we extend the notation for the characteristic functions
$\vfi,\vfit$ for conditions expressed in terms of $z_0,\z_0$.

\subsection{The density $\huu(x_0,y_0)$}\label{Shuu}

In this case $r=1$ and $\kk=k\ge 4$ is even, although the first part
of the calculation holds true more generally. Then, by
\eqref{eqCenter}, we know that 
$C_k:=C_\kk(1)=\big(x_0,\frac{x_0+y_0}{k}\big)$, 
and by the definition of
$\T_\k$ we know that 
$\C_k\in\overset{\circ}\T_k$ if and only if the following
conditions hold simultaneously:
    \begin{align*}
    \begin{cases}
     &1-\frac{x_0+y_0}{k}<x_0<1,\\
     & \frac{x_0+1}{k+1}<\frac{x_0+y_0}{k}<\frac{x_0+1}{k}.
     \end{cases}
    \end{align*}
Since we assumed that $x_0,y_0>0$, this translates into the equivalence
  \begin{equation*}
      C_k\in\overset{\circ}{\T_k}\ \ \iff\ \ 0<x_0, y_0<1
        \text{ and }     k<\frac{x_0+y_0}{1-\min(x_0,y_0)}\,,
          \quad\text{for $k\ge 2$.}
  \end{equation*}
For the edges of $\T_k$, we obtain
    \begin{align*}
        C_k(1)\in\partial{\T_k}\setminus V(\T_k)\ \ \iff\ \
    \begin{cases}
      k=\frac{x_0+y_0}{1-x_0},\quad 
                &\text{for\ \  $\frac{k-1}{k+1}< x_0<\frac{k}{k+2}$;} \\
       \qquad\text{or}& \\       
      k=\frac{x_0+y_0}{1-y_0},\quad 
                &\text{for\ \  $\frac{k}{k+2}< x_0<1$;} \\
       \qquad\text{or}& \\
      x_0=1,\quad &\text{for\ \  $\frac{k-1}{k+1}< y_0<1$;} \\
       \qquad\text{or}& \\       
      y_0=1,\quad &\text{for\ \  $\frac{k-1}{k+1}< x_0<1$,}
    \end{cases}
    \end{align*}
and for the vertices of $\T_k$, we have:
        \begin{equation*}
        C_k\in V(\T_k)\ \ \iff\ \  (x_0,y_0)\in
     \Big\{\Big(\frac{k-1}{k+1}, 1\Big);\  
                \Big(\frac k{k+2}, \frac k{k+2}\Big);\  
                  \Big(1, \frac{k-1}{k+1}\Big);\
                         (1, 1) 
                \Big\}=\U_{k}\,.
        \end{equation*}
Then, using Table~\ref{Table2}, for $\huu$, the corresponding sum from
the right-hand side of \eqref{eqgr} gives
     \begin{equation}\label{eqh1} 
        \begin{split}
     \huu(x_0,y_0) =&\phantom{+}
                \sum_{\substack{k\ge 4\\k\text{ even}}} 
        \frac 2k
                \vfit\Big(k<\frac{x_0+y_0}{1-z_0};\ 0<x_0,y_0<1\Big)\\
                &+\sum_{\substack{k\ge 4\\k\text{ even}}} 
        \frac 1{k}\bigg\{\vfi\Big(k=\frac{x_0+y_0}{1-z_0}, \text{ if }
                \frac{k-1}{k+1}<z_0<\frac{k}{k+2}\Big)\\
        &\phantom{\sum_{\substack{k\ge 4\\k\text{ even}}}\frac{1}{k}
                k=\frac{x_0+y_0}{1-z_0}\qquad}
        +\vfi\Big(z_0=1\ \text{ if }
                \frac{k-1}{k+1}<\z_0<1\Big)\bigg\}\\
        &+\sum_{\substack{k\ge 4\\k\text{ even}}} 
         \bigg\{\frac{2k+1}{4k(k+1)}
                \vfi\Big(z_0=\frac{k-1}{k+1};\ \z_0=1\Big)\\
        &\phantom{\sum_{\substack{k\ge 4\\k\text{ even}}}\frac{1}{k}
                k=\frac{x_0+y_0}{1-z_0}\qquad}
                +\frac{k+2}{2k(k+1)}
                \vfit\Big(z_0=\frac{k}{k+2}\Big)
                +\frac{2}{k}
                \vfit\Big(z_0=1\Big)
        \bigg\}\,.
        \end{split}
     \end{equation}
Here, for a given $(x_0,y_0)$, in each sum the number of nonzero terms is
finite, at most equal to one in the last two sums that correspond to
points on the border of $\U_k$.

\subsection{The density $\hdu(x_0,y_0)$}\label{Shdu}

Here $r=2$, $\kk=(1,l)$ or $\kk=(k,1)$,  with $k,l\ge 5$, both odd. 
We assume first that $\kk=(1,l)$. The center of
$\P_{1,l}(2)$ has coordinates 
$C_{1,l}:=C_{1,l}(2)=\big(x_0,\frac{lx_0+y_0}{l-1}\big)$, 
by \eqref{eqCenter}. Then
$\C_{1,l}\in\overset{\circ}\T_{1,l}$ if and only if the following
conditions hold simultaneously:
    \begin{align*}
    \begin{cases}
     &\frac{(l-1)\frac{lx_0+y_0}{l-1}-1}{l}<x_0<
        \frac{l\cdot\frac{lx_0+y_0}{l-1}-1}{l+1},\\
     & \frac{x_0+1}{2}<\frac{lx_0+y_0}{l-1}<1.
     \end{cases}
    \end{align*}
This gives
  \begin{equation*}
      C_{1,l}\in\overset{\circ}\T_{1,l}\ \ \iff\ \ 
    \begin{cases}
        &x_0<l,\ y_0<1,\\
        &\frac{1+y_0}{1-x_0}<l<\min\Big(
        \frac{x_0+y_0+1}{1-x_0},\;
        \frac{x_0+1}{1-y_0}\Big)\,,
     \end{cases}\quad\text{for $l\ge 3$.}
  \end{equation*}
Next, the conditions for the open edges of $\T_{1,l}$ are
    \begin{align*}
        C_{1,l}\in\partial{\T_{1,l}}\setminus V(\T_{1,l})\ \ \iff\ \
    \begin{cases}
      (l+1)x_0+2y_0=l-1,\quad 
                &\text{for\ \  $\frac{l-3}{l+1}< x_0<\frac{l-2}{l+2}$;} \\
       \qquad\text{or}& \\       
      x_0+ly_0=l-1,\quad 
                &\text{for\ \  $\frac{l-2}{l+2}< x_0<\frac{l-1}{l+1}$;} \\
       \qquad\text{or}& \\
     lx_0+y_0=l-1,\quad &\text{for\ \  $\frac{l-2}{l}< x_0<\frac{l-1}{l+1}$;} \\
       \qquad\text{or}& \\       
     y_0=1,\quad &\text{for\ \  $\frac{l-3}{l+1}< x_0<\frac{l-2}{l}$.}
    \end{cases}
    \end{align*}
For the vertices of $\T_{1,l}$, we have
        \begin{equation*}
        \begin{split}
        C_{1,l}\in V(\T_{1,l})\ \ \iff\ \  (x_0,y_0)\in
     \Big\{\Big(\frac{l-3}{l+1}, 1\Big);&\; 
                \Big(\frac{l-2}{l+2}, \frac{l}{l+2}\Big);\; \\
                &\Big(\frac{l-1}{l+1}, \frac{l-1}{l+1}\Big);\; 
                \Big(\frac{l-2}{l}, 1\Big)
        \Big\}=\U_{1,l}\,.
        \end{split}
        \end{equation*}
The symmetry allows us to use the $z_0,\z_0$ notation to collect the
contribution of terms corresponding to $(1,l)$ and $(k,1)$ in the same
formula. 
Then, using Table~\ref{Table2}, the corresponding sum from
the right-hand side of \eqref{eqgr} gives
     \begin{equation}\label{eqh2} 
        \begin{split}
     \hdu(x_0,y_0) =&\phantom{+}
                \sum_{\substack{l\ge 5\\k\text{ odd}}} 
        \frac 2{l-1}
                \vfi\Big(z_0,\z_0<1;\ 
        \frac{1+z_0}{1-\z_0}<l<
        \min\Big(\frac{2z_0+\z_0+1}{1-\z_0},\; \frac{\z_0+1}{1-z_0}\Big)\Big)\\
                &+\sum_{\substack{l\ge 5\\l\text{ odd}}} 
        \frac 1{(l-1)}\bigg\{\vfi\Big((l+1)z_0+2\z_0=l-1, \text{ if }
                \frac{l-3}{l+1}<z_0<\frac{l-2}{l+2}\Big)\\
        &\phantom{\sum_{\substack{k\ge 4\\k\text{ odd}}}\qquad\quad\qquad\quad}
        +\vfi\Big(z_0+l\z_0=l-1,\ \text{ if }
                \frac{l-2}{l+2}<\z_0<1\Big)\bigg\}\\
        &\phantom{\sum_{\substack{k\ge 4\\k\text{ odd}}}\qquad\quad\qquad\quad}
        +\vfi\Big(z_0=1,\ \text{ if }
                \frac{l-3}{l+1}<\z_0<\frac{l-2}{l}\Big)\bigg\}\\
        &+\sum_{\substack{l\ge 5\\l\text{ odd}}} 
         \bigg\{\frac{l}{2(l-1)(l+1)}
                \vfi\Big(z_0=\frac{l-3}{l+1};\ \z_0=1\Big)      \\
        &\phantom{\sum_{\substack{k\ge 5\\k\text{ odd}}}\frac{1}{2k}\quad}
                +\frac{2l^2+5l+1}{4(l-1)l(l-1)}
                \vfi\Big(z_0=\frac{l-2}{l+2};\ \z_0=\frac{l}{l+2}\Big)\\
        &\phantom{\sum_{\substack{k\ge 5\\k\text{ odd}}}\frac{1}{2k}\quad}
                +\frac{1}{l}
                \vfit\Big(z_0=\frac{l-1}{l+1}\Big)
        +\frac{2l+1}{4l(l-1)}
                \vfi\Big(z_0=\frac{l-2}{l};\ \z_0=1\Big)
        \bigg\}\,.
        \end{split}
     \end{equation}
Here, for a given $(x_0,y_0)$, in each sum the number of nonzero terms is
finite, at most equal to one in the last two sums that correspond to
points on the border of $\U_{1,l,1}$.

\subsection{The density $\htu(x_0,y_0)$}\label{Shtu}

Now $r=3$ and $\kk=(1,l,1)$, with $l\ge 6$ even. The center of
$\P_{1,l,1}(3)$ is 
$C_{1,l,1}:=C_{1,l,1}(3)=\big(x_0,\frac{(l-1)x_0+y_0}{l-1}\big)$, 
cf. \eqref{eqCenter}, and $\T_{1,l,1}=\T_{1,l}$. Then
$\C_{1,l,1}\in\overset{\circ}\T_{1,l,1}$ if and only if the following
conditions hold simultaneously:
    \begin{align*}
    \begin{cases}
     &\frac{(l-1)\frac{(l-1)x_0+y_0}{l-2}-1}{l}<x_0<
        \frac{l\cdot\frac{(l-1)x_0+y_0}{l-2}-1}{l+1},\\
     & \frac{x_0+1}{2}<\frac{(l-1)x_0+y_0}{l-2}<1.
     \end{cases}
    \end{align*}
This can be rewritten as
  \begin{equation*}
      C_{1,l,1}\in\overset{\circ}\T_{1,l,1}\ \ \iff\ \ 
    \begin{cases}
        x_0-y_0+2&<l(1-y_0),\\
        l(1-y_0)&<2(1+x_0),\\
        l(1-x_0)&<2(y_0+1),\\
        y_0-x_0+2&<l(1-x_0).    
     \end{cases}
  \end{equation*}
The conditions for the open edges of $\T_{1,l,1}$ are
    \begin{align*}
        C_{1,l,1}\in\partial{\T_{1,l,1}}\setminus V(\T_{1,l,1})\ \ \iff\ \
    \begin{cases}
      lx_0+2y_0=l-2,\quad 
                &\text{for\ \  $\frac{l-3}{l+1}< x_0<\frac{l-2}{l+2}$;} \\
       \qquad\text{or}& \\       
      2x_0+ly_0=l-2,\quad 
                &\text{for\ \  $\frac{l-2}{l+2}< x_0<\frac{l-1}{l+1}$;} \\
       \qquad\text{or}& \\
   (l-1)x_0+y_0=l-2,\quad &\text{for\ \  $\frac{l-2}{l}< x_0<\frac{l-1}{l+1}$;} \\
       \qquad\text{or}& \\       
     x_0+(l-1)y_0=l-2,\quad &\text{for\ \  $\frac{l-3}{l+1}< x_0<\frac{l-2}{l}$.}
    \end{cases}
    \end{align*}
For the vertices of $\T_{1,l,1}$, we have:
        \begin{equation*}
        \begin{split}
        C_{1,l,1}\in V(\T_{1,l,1})\ \ \iff\ \  (x_0,y_0)\in
            \Big\{\Big(\frac{l-3}{l+1},& \frac{l-1}{l+1}\Big);\; 
                \Big(\frac{l-2}{l+2}, \frac{l-2}{l+2}\Big);\;\\ 
                &\Big(\frac{l-1}{l+1}, \frac{l-3}{l+1}\Big);\; 
                \Big(\frac{l-2}{l}, \frac{l-2}{l}\Big)
        \Big\}=\U_{1,l,1}\,.
        \end{split}
        \end{equation*}
Then, using Table~\ref{Table2}, the corresponding sum from
the right-hand side of \eqref{eqgr} gives
     \begin{equation}\label{eqh3} 
        \begin{split}
     \htu(x_0,y_0) =&\phantom{+}
                \sum_{\substack{l\ge 6\\k\text{ even}}} 
        \frac 2{l-2}
                \vfit\Big(z_0-\z_0+2<l(1-\z_0);\ l(1-z_0)<2(\z_0+1)\Big)\\
                &+\sum_{\substack{l\ge 6\\l\text{ even}}} 
        \frac 1{(l-2)}\bigg\{\vfi\Big(lz_0+2\z_0=l-2, \text{ if }
                \frac{l-2}{l+2}<z_0<\frac{l-3}{l+1}\Big)\\
        &\phantom{\sum_{\substack{k\ge 4\\k\text{ even}}}\qquad\quad\qquad\quad}
        +\vfi\Big(z_0+(l-1)\z_0=l-2,\ \text{ if }
                \frac{l-3}{l+1}<\z_0<\frac{l-2}{l}\Big)\bigg\}\\
        &+\sum_{\substack{l\ge 6\\l\text{ even}}} 
         \bigg\{\frac{2l-1}{4l(l-1)}
                \vfi\Big(z_0=\frac{l-3}{l+1};\ \z_0=\frac{l-1}{l+1}\Big)\\
        &\phantom{\sum_{\substack{k\ge 4\\k\text{ even}}}\frac{1}{2k}\quad}
                +\frac{l+2}{2l(l-2)}
                \vfit\Big(z_0=\frac{l-2}{l+2}\Big)
                +\frac{l}{2(l-2)(l-1)}
                \vfit\Big(z_0=\frac{l-2}{l}\Big)
        \bigg\}\,.
        \end{split}
     \end{equation}
Here, for a given $(x_0,y_0)$, in each sum the number of nonzero terms is
finite, at most equal to one in the last two sums that correspond to
points on the border of $\U_{1,l,1}$.

\subsection{The contribution to the density of points of type
$\TT(r)$, $r\ge5$}\label{Sgu}
In this section we find the tail density, which we define to be
        \begin{equation*}
                g^\mathrm{u}(x_0,y_0)=\sum_{r\ge5}g_r(x_0,y_0)\,.
        \end{equation*} 

For any $r\ge5$, there exist only two admissible $r-$tuples $\kk$,
symmetric to one another. Let $\kk=(1,2,\dots,2,3)$ be such an
$r-$tuple. Then, since $p_r(\kk)=2$ and $p_{r-1}(2,\dots,3)=2r-1$, 
relation \eqref{eqCenter} 
produces $C_\kk:=C_\kk(3)=\big(x_0,\frac{(2r-1)x_0+y_0}{2}\big)$.
Then $\C_\kk\in\overset{\circ}\T_\kk$ if and only if the following
conditions hold simultaneously:
    \begin{align*}
    \begin{cases}
     &\frac{2\frac{(2r-1)x_0+y_0}{2}-1}{2r-1}<x_0<
        \frac{2\frac{(2r-1)x_0+y_0}{2}-1}{2r-3},\\
     & 1-x_0<\frac{(2r-1)x_0+y_0}{2}<1.
     \end{cases}
    \end{align*}
This can be rewritten as
  \begin{equation*}
      C_\kk\in\overset{\circ}\T_\kk\ \ \iff\ \ 
    \begin{cases}
        &\!\!\!y_0<1,\\
        &1<2x_0+y_0,\\
        &2<(2r+1)x_0+y_0,\\
        &2>(2r-1)x_0+y_0.       
     \end{cases}
  \end{equation*}
The conditions for the open edges of $\T_\kk$ are
    \begin{align*}
        C_\kk\in\partial{\T_\kk}\setminus V(\T_\kk)\ \ \iff\ \
    \begin{cases}
      (2r+1)x_0+y_0=2,\quad 
                &\text{for\ \  $\frac{1}{2r+1}< x_0<\frac{2r-1}{l+2}$;} \\
       \qquad\text{or}& \\       
      2x_0+y_0=1,\quad 
                &\text{for\ \  $\frac{1}{2r-1}< x_0<\frac{1}{2r-3}$;} \\
       \qquad\text{or}& \\
   (2r-1)x_0+y_0=2,\quad &\text{for\ \  $\frac{1}{2r-1}< x_0<\frac{1}{2r-3}$;} \\
       \qquad\text{or}& \\       
     y=1,\quad &\text{for\ \  $\frac{1}{2r+1}< x_0<\frac{1}{2r-1}$.}
    \end{cases}
    \end{align*}
For the vertices of $\T_\kk$, we have
        \begin{equation*}
        \begin{split}
        C_\kk\in V(\T_\kk)\ \ \iff\ \  (x_0,y_0)\in
        \Big\{\Big(\frac{1}{2r+1}, 1\Big);\;& 
                \Big(\frac{1}{2r-1}, \frac{2r-3}{2r-1}\Big);\; \\
                &\Big(\frac{1}{2r-3}, \frac{2r-5}{2r-3}\Big);\; 
                \Big(\frac{1}{2r-1}, 1\Big)
        \Big\}=\U_\kk\,.
        \end{split}
        \end{equation*}
We use the notation with the variables $z_0,\z_0$ and the symmetry to
write in the same formula the contribution to $g(x_0,y_0)$ of all the
terms corresponding to the $r-$tuples $(1,2,\dots,2,3)$ and
$(3,2,\dots,2,1)$, for $r\ge5$. Using  the information from 
Table~\ref{Table2}, we obtain
     \begin{equation}\label{eqhr} 
        \begin{split}
   g^\mathrm{u}(x_0,y_0) =&\phantom{+}
                \sum_{r\ge 5}
                \vfi\big(\z_0<1<2z_0+\z_0;\ (2r-1)z_0+\z_0<2<(2r+1)z_0+\z_0\big)\\
                &+\sum_{r\ge 5} \frac 12\bigg\{
                   \vfi\Big((2r+1)z_0+\z_0=2, \text{ if }
                \frac{1}{2r+1}<z_0<\frac{1}{2r-1}\Big)\\
        &\qquad\qquad\qquad\qquad+
                   \vfi\Big(2z_0+\z_0=1, \text{ if }
                \frac{1}{2r-1}<z_0<\frac{1}{2r-3}\Big)\\
        &\qquad\quad+
                        \vfi\Big((2r-1)z_0+\z_0=2, \text{ if }
                \frac{1}{2r-1}<z_0<\frac{1}{2r-3}\Big)\\
        &\qquad\qquad\qquad\qquad+
                        \vfi\Big(\z_0=1, \text{ if }
                \frac{1}{2r+1}<z_0<\frac{1}{2r-1}\Big)
                                \bigg\}\\
        &+\sum_{r\ge 5} 
         \bigg\{\frac{4r+1}{8(2r+1)}
                \vfi\Big(z_0=\frac{1}{2r+1};\ \z_0=1\Big)\\
        &\qquad\qquad\qquad\qquad
                +\frac{14r+9}{16(2r+1)}
                \vfi\Big(z_0=\frac{1}{2r-1};\ \z_0=\frac{2r-3}{2r-1}\Big)\\
        &\qquad\quad+
                \frac{2r-3}{16(2r-1)}
                \vfi\Big(z_0=\frac{1}{2r-3};\ \z_0=\frac{2r-5}{2r-3}\Big)\\
        &\qquad\qquad\qquad\qquad
                +\frac{4r-1}{8(2r-1)}
                \vfi\Big(z_0=\frac{1}{2r-1};\ \z_0=1\Big)
        \bigg\}\,.
        \end{split}
     \end{equation}
Here, in each sum, for a given $(x_0,y_0)$, at most one term is nonzero.

\subsection{The baby puzzle}\label{Sbaby}
Let $i\ge 5$ be a positive integer. Then, we observe that the polygons
$\U_{i-1}$, $\U_{1,i}$, $\U_{i,1}$, $\U_{1,i,1}$ fit perfectly into
the quadrilateral with vertices
        \begin{equation*}
                \W_i=\Big\{
                        \Big(\frac{i-3}{i+1}, 1\Big);\  
                \Big(\frac {i-1}{i+3}, \frac {i-1}{i+3}\Big);\  
                  \Big(1, \frac{i-3}{i+1}\Big);\ (1, 1) 
                \Big\}\,,
        \end{equation*}
(see Figure~\ref{Figure6}).
Another nice aspect of this matching is due to the
fact that $p_1(i-1)=p_2(1,i)=p_2(i,1)=p_3(1,i+1,1)$, so $\W_i$ may be
viewed as region with constant density at level $i$, say. 
Let us notice that $\W_5\supset\W_7\supset\W_9\supset,\dots$, so
$g(u,v)$ is not locally constant on $\W_i$, but the support of
$g(u,v)$ is a superposition of quadrilateral steps of constant
density (provided we show that a similar property holds for the remaining
polygons $\U_\kk$). 
Putting together \eqref{eqh1}, \eqref{eqh2}, \eqref{eqh3}, we obtain
the next result.

\begin{proposition}\label{Proposition1}
For any $(x_0,y_0)\in[0,1]^2$, we have
     \begin{equation}\label{eqP1} 
        \begin{split}
    h^\mathrm{u}(x_0,y_0):=& \huu(x_0,y_0)+\hdu(x_0,y_0)+\htu(x_0,y_0) =\\
        =&\phantom{+}
                \sum_{\substack{i\ge 5\\i\mathrm{\ odd}}} 
        \frac 2{i-1}
                \vfit\Big(z_0<1;\ i<\frac{z_0+2\z_0+1}{1-z_0}\Big)\\
                &+\sum_{\substack{i\ge 5\\i\mathrm{\ odd}}} 
        \frac 1{(i-1)}\bigg\{\vfi\Big((i+1)z_0+2\z_0=i-1, \mathrm{\ if\ }
                \frac{i-3}{i+1}<z_0<\frac{i-1}{i+3}\Big)\\
        &\phantom{\sum_{\substack{k\ge 4\\k\text{ even}}}\frac{1}{2k}
                k=\frac{x_0+y_0}{1-z_0}\qquad\qquad}
        +\vfi\Big(z_0=1,\ \mathrm{\ if\ }
                \frac{i-3}{i+1}<\z_0<1\Big)\bigg\}\\
        &+\sum_{\substack{i\ge 5\\i\mathrm{\ odd}}} 
         \bigg\{\frac{i}{2(i-1)(i+1)}
                \vfi\Big(z_0=\frac{i-3}{i+1};\ \z_0=1\Big)\\
        &\phantom{\sum_{\substack{k\ge 4\\k\text{ even}}}\frac{1}{2k}++}
                +\frac{i+3}{2(i-1)(i+1)}
                \vfit\Big(z_0=\frac{i-1}{i+3}\Big)
                +\frac{1}{2(i-1)}
                \vfit\Big(z_0=1\Big)
        \bigg\}\,.
        \end{split}
     \end{equation}
\end{proposition}

\begin{proof}
We only need to check the equality at the matching corners. We have:
  \begin{align*}
        h^\mathrm{u}\bigg(\frac{i-2}{i+2},\frac{i}{i+2}\bigg)
        =&\frac 12\Big(\alpha_{1,i}\big(\tfrac{i-2}{i+2},\tfrac{i}{i+2}\big)+
        \alpha_{1,i+1,1}\big(\tfrac{i-2}{i+2},\tfrac{i}{i+2}\big)\Big)\\
        =&\frac 12\bigg(\frac{2i^2+5i+1}{2(i-1)i(i+1)}+\frac{2i+1}{2i(i+1)}\bigg)
        =\frac{1}{i-1}\,,\\
\intertext{and}
        h^\mathrm{u}\bigg(\frac{i-2}{i},1\bigg)
        =&\frac 12\Big(\alpha_{1,i}\big(\tfrac{i-2}{i},1\big)+
        \alpha_{i-1}\big(\tfrac{i-2}{i},\tfrac{2}{i}\big)\Big)\\
        =&\frac 12\bigg(\frac{2i+1}{2i(i-1)}+\frac{2i-1}{2i(i-1)}\bigg)
        =\frac{1}{i-1}\,,
  \end{align*}
equal each to half of the interior density, since they are on the open
edges of $\W_i$. By symmetry, we have the same result at 
$\big(i/(i+2),(i-2)/(i+2)\big)$ and at 
$\big(1,(i-2)/i\big)$. At the interior matching point, we have
  \begin{equation*}
    \begin{split}
        h^\mathrm{u}\bigg(\frac{i-1}{i+1},\frac{i-1}{i+1}\bigg)
        =&\frac 12\Big(\alpha_{i-1}\big(\tfrac{i-1}{i+2},\tfrac{2}{i+1}\big)+
        \alpha_{1,i}\big(\tfrac{i-1}{i+1},1\big)+
        \alpha_{i,1}\big(\tfrac{i-1}{i+1},\tfrac{2}{i+1}\big)+
        \alpha_{1,i+1,1}\big(\tfrac{i-1}{i+1},1\big)\Big)\\
        =&\frac 12\bigg(\frac{i+1}{i(i-1)}+\frac{2}{i}+\frac{i+1}{i(i-1)}\bigg)
        =\frac{2}{i-1}\,,
    \end{split}
  \end{equation*}
which concludes the proof of the proposition.
\end{proof}

\subsection{The big puzzle}\label{Sbig}
We group the terms of lower orders into
  \begin{equation*}
    \begin{split}
        g^\mathrm{d}(x_0,y_0)=h^\mathrm{d}(x_0,y_0)+h^\mathrm{u}(x_0,y_0)+
                        g_4(x_0,y_0),
    \end{split}
  \end{equation*}
where
$h^\mathrm{d}(x_0,y_0)=\hud(x_0,y_0)+\hdd(x_0,y_0)+\htd(x_0,y_0)$. It
remains to find $h^\mathrm{d}(x_0,y_0)+ g_4(x_0,y_0),$ as
$h^\mathrm{u}(x_0,y_0)$ was the object of Section~\ref{Sbaby}. In
fact, the calculations are special cases of those already performed in
sections~\ref{Shuu}-\ref{Sgu}. Here, we shall show that the sum of
$h^\mathrm{d}$, $g_4$ and $g^\mathrm{u}$ can be combined into a
simpler formula, similar to that of $h^\mathrm{u}$. 

Let 
$\M=\big\{(2)$;
$(1,3)$; $(3,1)$;
$(1,2,3)$; $(3,2,1)$; 
$(1,4,1)$; $(1,2,4,1)$; $(1,4,2,1)$; 
$(1,2,2,3)$; $(3,2,2,1)$;
$(1,2,2,2,3)$; $(3,2,2,2,1);\dots\big\}$.
The key point in the matching that occurs among the supports of
different components of $h^\mathrm{d}+g_4+g^\mathrm{u}$ is the fact
that the kernel is constant, equal to $1$, for any $\kk$
contributing to the sum. This follows by the equality $p_r(\kk)=2$,
for any $\kk\in\M$.

Let $\W$ be the quadrilateral with vertices 
$(0,1);$ $(1/3,1/3);$ $(1,0);$ $(1,1)$. It turns out that $\W$ is the
support of $h^\mathrm{d}+g_4+g^\mathrm{u}$, and it looks like a mosaic (that
is, no superpositions over interior points occur) composed by all polygons
$\U_\kk$, with $\kk\in\M$. Perfect matchings occur (see
Figure~\ref{Figure7}) getting particular cases
of the general relations, as follows. (In the nonsymmetric cases we
give the statement only for the polygons situated above the first
diagonal.) 
The polygons $U_{1,2,2,3}$ and $U_{1,2,3}\cup\U_{1,2,4,1}$ are 
given also by the formula for $\U_{1,2,\dots,2,3}$, with $r=4$ and
$r=3$, respectively. More significantly, on the one hand
$\wU_{1,4,1}=\U_{1,2,4,1}\cup\U_{1,4,1}\cup\U_{1,4,2,1}$ is given by
the formula for $\U_{1,l,1}$, with $l=4$ and, on the other hand, 
$\wU_{1,3}=\U_{1,3}\cup\U_{1,2,3}\cup\U_{1,2,2,3}\cup\U_{1,2,2,2,3}\cdots$
is given by the formula for $\U_{1,l}$, with $l=3$. 
Then $\W$ is composed by $\wU_{1,3}$, $\wU_{1,4,1}$, $\wU_{3,1}$, and
$\U_{1}$ in the same way as $\U_{1,i}$, $\U_{1,i+1,1}$, $\wU_{i,1}$, and
$\U_{i-1}$ completed the baby puzzle, $\W_i$. Next we summarize the
contribution to $g(x_0,y_0)$ of all tuples $\kk\in\M$. 

\begin{proposition}\label{Proposition2}
For any $(x_0,y_0)\in[0,1]^2$, we have
     \begin{equation}\label{eqP2} 
        \begin{split}
        g^\mathrm{du}(x_0,y_0):&=
        h^\mathrm{d}(x_0,y_0)+g_4(x_0,y_0)+g^\mathrm{u}(x_0,y_0)=\\
                &=\vfit\big(z_0<1;\ 1<2z_0+\z_0\big)+\\
        &+\frac 12 \vfi\big(2z_0+\z_0=1, \mathrm{\ if\ }
                0<z_0<1/3\big)
        +\frac 12\vfi\big(z_0=1, \mathrm{\ if\ }
                0<z_0<1\big)\\
        &+\frac{3}{16}\vfi(z_0=0;\ \z_0=1)
        +\frac{3}{8}\vfit(z_0=1/3)
        +\frac{1}{4}\vfit(z_0=1)\,.
        \end{split}
     \end{equation}
\end{proposition}

\begin{proof}
We just check the equality at the matching vertices using entries from
Table~\ref{Table2}. Beginning with the border of $\W$,
for points on the top edge, we have:
  \begin{align*}
        g^\mathrm{du}(1/3,1)=&\frac 12\big(\alpha_{2}(1/3,2/3)+
                        \alpha_{13}(1/3,1)\big)
                =\frac 12\bigg(\frac{5}{12}+\frac{7}{12}\bigg)
                        =\frac 12\,;\\
        g^\mathrm{du}(1/5,1)=&\frac 12\big(\alpha_{13}(1/5,4/5)+
                        \alpha_{123}(1/5,1)\big)
                =\frac 12\bigg(\frac{9}{20}+\frac{11}{20}\bigg)
                        =\frac 12\,;\\
        g^\mathrm{du}(1/7,1)=&\frac 12\big(\alpha_{123}(1/7,6/7)+
                        \alpha_{1223}(1/7,1)\big)
                =\frac 12\bigg(\frac{13}{28}+\frac{15}{28}\bigg)
                        =\frac 12\,;\\
        g^\mathrm{du}\big(\tfrac{1}{2r+1},1\big)
                =&\frac 12\Big(
                \alpha_{12\cdots23}\big(\tfrac{1}{2r+1},\tfrac{2r}{2r+1}\big)+
                        \alpha_{122\cdots23}\big(\tfrac{1}{2r+1},1\big)\Big)\\
                =&\frac 12\bigg(
                        \frac{4r+1}{4(2r+1)}+   \frac{4(r+1)-1}{4(2(r+1)-1)}
                \bigg)=\frac 12,\,\quad\text{for $r\ge 4$}.
   \end{align*}
On the edge with endpoints $(0,1)$, $(1/3,1/3)$, we have:
  \begin{align*}
        g^\mathrm{du}(1/5,3/5)=&\frac 12\big(\alpha_{123}(1/5,4/5)+
                                \alpha_{1223}(1/5,1)+
                                \alpha_{1241}(1/5,4/5)\big)\\
                =&\frac 12\bigg(
                        \frac{13}{21}+
                        \frac{5}{56}+\frac{7}{24}
                \bigg)=\frac 12\,;\\
        g^\mathrm{du}(1/7,5/7)=&\frac 12\big(\alpha_{1223}(1/7,6/7)+
                                \alpha_{12223}(1/7,1)\big)
                =\frac 12\bigg(
                        \frac{65}{72}+
                        \frac{7}{72}
                \bigg)=\frac 12\,;\\
        g^\mathrm{du}\big(\tfrac{1}{2r-1},\tfrac{2r-3}{2r-1}\big)
                =&\frac 12\Big(
                \alpha_{12\cdots23}\big(\tfrac{1}{2r-1},\tfrac{2(r-1)}{2r-1}\big)+
                        \alpha_{122\cdots23}\big(\tfrac{1}{2r-1},1\big)\bigg)\\
                =&\frac 12\bigg(
                        \frac{14r+9}{8(2r+1)}+  \frac{2(r+1)-3}{8(2(r+1)-1)}
                \Big)=\frac 12,\,\quad\text{for $r\ge 5$}.
   \end{align*}
These are all equal to half of the kernel, as expected.
The same results hold for points symmetric with respect
to the first diagonal.

At the vertices of $\W$, we have:
  \begin{align*}
        g^\mathrm{du}(1,1)=&\frac 12\alpha_{2}(1,1)=\frac 12\cdot\frac
                12=\frac 14\,;\\
        g^\mathrm{du}(1/3,1/3)=&\frac 12\big(\alpha_{141}(1/3,2/3)+
                                \alpha_{1241}(1/3,1)+
                                \alpha_{1421}(1/3,2/3)\big)\\
                =&\frac 12\bigg(
                         \frac{3}{5}+\frac{3}{40}+\frac{3}{40}
                \bigg)=\frac 3{18}\,.
   \end{align*}
At the corners $(0,1)$ and $(1,0)$ the problem is a little bit more
complicated. Here, as the inside
limit in \eqref{eqdefgr} is taken over $Q$, the parallelogram
$\P_\kk(r)$ will cover completely $\T_\kk$, for infinitely many $\kk$
from the sequences with $2$'s embraced by $1$ and $3$, while some
of these $\T_\kk$'s will be covered partially. 
Though,  we have a quick shot solution to
the problem of finding $g^\mathrm{du}(0,1)$ and $g^\mathrm{du}(1,0)$ 
using, on a larger scale, the property
of the parallelogram used in \eqref{eqremarq} for the quadrilateral with
vertices: $(0,1)$; $(1/3,1/3)$; $(1/2,1/2)$; $(1/3,1)$ and the fact
that in $\W$ the index has everywhere the same value $2$. We obtain
   \begin{equation*}
        \begin{split}
        g(0,1)=&\frac 12\bigg(2-\frac 12\alpha_{141}(1/2,1)
                        -\alpha_{13}(1/2,1)-\alpha_{13}(1/3,1)\\
                &\phantom{\frac 14\bigg(2-\frac 12\alpha_{141}(1/2,1/2)}
                        -\alpha_{1241}(1/3,1)
                                -\frac 12\alpha_{141}(1/3,2/3)\bigg)\\
        =&\frac 12\bigg(2-\frac 12\cdot \frac{2}{3}-\frac 13-\frac 7{12}
                -\frac 3{40}-\frac 12\cdot\frac 35\bigg)
                =\frac 3{16}\,.
        \end{split}
   \end{equation*}

Finally, we complete the puzzle with $\T_{1,4,1}$, the last piece. We have
  \begin{align*}
        g^\mathrm{du}(1/2,1/2)=&
        \frac 12\big(\alpha_{2}(1/2,1/2)+\alpha_{13}(1/2,1)
                +\alpha_{31}(1/2,1/2)+\alpha_{141}(1/2,1)\big)\\
                =&\frac 12\bigg(\frac 23+\frac 13+\frac 13+ \frac 23
                        \bigg)=1,\\
        g^\mathrm{du}(2/7,4/7)=&\frac 12\big(\alpha_{13}(2/7,5/7)+
                                \alpha_{123}(2/7,1)+
                                \alpha_{141}(2/7,5/7)+
                                \alpha_{1241}(2/7,1)\big)\\
                =&\frac 12\bigg(
                        \frac{19}{30}+
                         \frac{11}{30}+\frac{11}{30}+
                        \frac{19}{30}
                \bigg)=1,
   \end{align*}
and, by symmetry, $g^\mathrm{du}(4/7,2/7)=1$. We see that
$g^\mathrm{du}$ takes at $(1/2,1/2)$, $(2/7,4/7)$ and $(4/7,2/7)$ the
value $1$, equal to the kernel, as needed, since these are interior points.
This completes the proof of the proposition.

\end{proof}

\subsection{Completion of the proof of Theorem~\ref{Theorem1}}
We first notice that $\W$ can be obtained as a particular case of the
formula for $\W_i$, with $i=3$. Moreover, the right-hand side of
\eqref{eqP2} can also be obtained if we put $i=3$ into the generic terms
of the sums on the right-hand side of \eqref{eqP1}. Thus, we employ
Propositions~\ref{Proposition1} and~\ref{Proposition2} to obtain the sum 
$g(x_0,y_0)=g^\mathrm{du}(x_0,y_0)+h^\mathrm{u}(x_0,y_0)$, which
completes the proof of Theorem~\ref{Theorem1}.

\section{Proof of Theorem~\ref{TheoremI}}
Let $\I\subseteq [0,1]$ be fixed and let $a'/q',a''/q''$ be consecutive
fractions in $\FQ$. 
Let us first see how one translates the condition $a'/q'\in\I$ in
terms of the two variables $q',q''$. Since $a''q'-a'q''=1$, it follows
that $a'\equiv -\overline{q''}\pmod{q'}$. Here $\overline{q''}$ is
the representative from $[0,q'-1]$ of the inverse of $q''$ modulo
$q'$. Then, one immediately derives that $a'/q'\in\I$ if and only if
$q'-\overline{q''}\in q'\I$. Similarly, we get that $a''/q''\in\I$ if
and only if $\overline{q'}\in q''\I$, but here the inverse is taken
modulo $q''$. We remark that one of these conditions is almost redundant, since
$a'/q'\in\I$ ensures that $a''/q''\in\I$, also,  except for
at most one pair $q',q''$, and conversely. Then, imposing only one of
these two conditions, one may neglect this at most one term in the
corresponding computations below, and absorb it in the error term.  

The proof of Theorem~\ref{TheoremI} follows the same steps from the
beginning of the proof of Theorem~\ref{Theorem1}. We have to find the
ratio of the number of elements in the set
  \begin{equation*}
     \B_Q^\I(r)=\left\{(q',q'')\in \NN^2\colon\ 
        \begin{array}{l}  1\le q',q''\le Q,\ \gcd(q',q'') = 1,\ q'+q''>Q,\
                        \overline{q'}\bmod q''\in q''\I
                 \\ \displaystyle
        q'~ \text{even},\ q''~ \text{odd}; \ \kk(q',q'')\in\A(r),\ 
              (q',q^\LL(r))\in Q\cdot \Box
        \end{array} \right\},
  \end{equation*}  
and the cardinality of $\FQeI$. Then the turning point
is the analogue of relation \eqref{eqRaportul}, which becomes
   \begin{equation}\label{eqRaportulI}
        \iint\limits_{\Box_\eta(x_0,y_0)}g_r^\I(x,y)\,dxdy
         =\lim_{Q\rightarrow \infty}
           \frac{\#\B_Q^\I(r)}{\#\FQeI-1}\,.
   \end{equation}
This will complete the proof, provided we show that 
$\#\B_Q^\I(r)\sim |\I|\cdot\#\B_Q(r)$ and $\FQeI\sim |\I|\cdot\#\FQe$,
as $Q\to\infty$. One should observe that, with the notations from
Section~\ref{secProofT}, we have $\B_Q(r)=\B_Q^{[0,1]}(r)$.

To proceed, we estimate $\FQeI$.

\begin{lemma}\label{LemmaFQeI}
For any subinterval $\I\subseteq [0,1]$, we have
  \begin{equation}\label{eqFQeI}
        \begin{split}
     \FQeI&=\frac{|\I|Q^2}{\pi^2}+O\left(Q^{3/2+\varepsilon}\right)\,,
        \end{split}
  \end{equation}  
\end{lemma}

\begin{proof}

The cardinality of $\FQeI$ can be written as
  \begin{equation}\label{eqFeI1}
        \begin{split}
     \#\FQeI&=\#\left\{(q',q'')\in \NN^2\colon\ 
        \begin{array}{l}  1\le q',q''\le Q,\ \gcd(q',q'') = 1,\ q'+q''>Q,
                 \\ \displaystyle
        q'~ \text{even},\ q''~ \text{odd},\ \overline{q'}\bmod q''\in q''\I
        \end{array} \right\}\\
        &=\sum_{\substack{1\le q\le Q\\ q\text{ odd}}}\#\big\{
                x\in(Q-q,Q]\colon \ \gcd(x,q)=1,\ x\ \text{even, }\overline{x}\bmod
                q\in q\I\big\}\,.
        \end{split}
  \end{equation}  
Notice that the condition $\overline{q'}\in q''\I$ introduces
randomness in the positioning of points from the right hand side of
\eqref{eqFeI1}. In order to estimate the terms added in the sum, we
write them using exponential sums, separate the main term, and employ
classical bounds for Kloosterman sums (see \cite{Esterman},
\cite{Weil}). Thus, they can be written as 
  \begin{equation*}
        \begin{split}
        \sum_{\substack{Q-q<x\le Q\\ \gcd(x,q)=1\\ x\text{ even}}}\ \ &
        \sum_{y\in q\I}
        \frac 1q\sum_{k=1}^qe\Big(k\frac{y-\overline x}{q}\Big)=\\
                &=\frac{\phi(q)}{q}\cdot \frac{q}{2}\cdot q|\I|
                        +\frac 1q\sum_{k=1}^{q-1}
                                \sum_{y\in q\I} e\Big(k\frac{y}{q}\Big)         
                \sum_{\substack{Q-q<x\le Q\\ \gcd(x,q)=1\\ x\text{ even}}}
                        e\Big(\frac{-k\overline x}{q}\Big)\\
                &=\frac{\phi(q)}{q}\cdot
                        \frac{q}{2}\cdot q|\I|+O\left(q^{1/2+\varepsilon}\right)\,,
        \end{split}
  \end{equation*}  
in which $\phi(q)$ is Euler's totient function. Then, substituting in
\eqref{eqFeI1}, we obtain
  \begin{equation}\label{eqFeI2}
        \begin{split}
     \FQeI&=\frac{|\I|}{2}
        \sum_{\substack{1\le q \le Q\\ q\text{ odd}}}
                \frac{\phi(q)}{q}\cdot q+O\left(Q^{3/2+\varepsilon}\right)\,,
        \end{split}
  \end{equation}  
It remains to estimate the sum from \eqref{eqFeI2}. This is
  \begin{equation}\label{eqFeI3}
        \begin{split}
        \sum_{\substack{1\le q \le Q\\ q\text{ odd}}}\frac{\phi(q)}{q}\cdot q
        &=\sum_{\substack{1\le q \le Q\\ q\text{ odd}}}
                \sum_{d\mid q}\frac{\mu(d)}{d}\cdot q
        =\sum_{\substack{d=1\\ d\text{ odd}}}^Q\frac{\mu(d)}{d}
        \sum_{\substack{q_1=1\\ q_1\text{ odd}}}^{Q/d}d\cdot q_1\\
        &=\sum_{\substack{d=1\\ d\text{ odd}}}^Q\frac{\mu(d)}{d}
                \cdot\frac{1}{2d}\int_1^Q t\,dt+O(Q\log Q)\\
                &=\frac{2Q^2}{\pi^2}+O(Q\log Q),
        \end{split}
  \end{equation}  
since, via the Euler product, we find that  
  \begin{equation*}
        \begin{split}
        \sum_{\substack{d=1\\ d\text{ odd}}}^Q\frac{\mu(d)}{d^2}
        &=\prod_{p\ge 3}\left(1-\frac 1{p^2}\right)+O(Q)\\
        &=\frac{6}{\pi^2}\left(1-\frac 1{2^2}\right)^{-1}+O(Q).
        \end{split}
  \end{equation*}  
Now, the lemma follows by inserting the estimation \eqref{eqFeI3} on the right-hand side of
\eqref{eqFeI2}.
\end{proof}

We remark that the size of the error term in \eqref{eqFQeI} may be slightly lowered,
but this is not essential for our needs.

Next, for any $\Omega\subset \RR^2$ and $\I\subseteq [0,1]$, we denote 
   \begin{equation*}
          \Neo^\I(\Omega):=\#\left\{(x,y)\in \Omega\cap
              \ZZ^2\colon~x~\text{even},~y~\text{odd},\ \gcd(x,y) = 1,\
                        \overline{x}\bmod y \in y\I
         \right\}.
   \end{equation*}

Then, in the spirit of Lemma~\ref{LemmaFQeI}, we get, more generally, the following result.

\begin{lemma}\label{LemmaNeoI}
Let $R>0$ and $\Omega\subseteq [0,R]\times [0,R]$ be a convex domain.
Then, we have
   \begin{equation*}
  \Neo^\I(\Omega)= |\I|\cdot \Neo(\Omega) + O(R^{3/2+\varepsilon})\,.
   \end{equation*}
\end{lemma}
\begin{proof}
The proof follows the first part of the proof of Lemma~\ref{LemmaFQeI}.
We have:
  \begin{equation}\label{eqNFeI1}
        \begin{split}
     \Neo^\I(\Omega)&=\#\left\{(q_1,q_2)\in\Omega\cap\NN^2\colon\ 
        \begin{array}{l}  \gcd(q_1,q_2) = 1,\ q_1~ \text{even},\ q_2~ \text{odd},
                 \\ \displaystyle
                        \overline{q_1}\bmod q_2\in q_2\I
        \end{array} \right\}\\
        &=\sum_{\substack{1\le q\le R\\ q\text{ odd}}}\#\big\{
                x\in \mathfrak{I}(q)\colon \ \gcd(x,q)=1,\ x\ \text{even, }\overline{x}\bmod
                q\in q\I\big\}\,,
        \end{split}
  \end{equation}  
in which $\mathfrak{I}(q):=\Omega\cap\{ y=q\}$.
Evaluating the size of the terms here, using again the same estimates for
Kloosterman sums, we find that they are equal to
  \begin{equation}\label{eqNFeI2}
        \begin{split}
        \sum_{\substack{x\in\mathfrak{I}(\Omega)\\ \gcd(x,q)=1\\ x\text{ even}}}\ \
        \sum_{y\in q\I}
        \frac 1q\sum_{k=1}^qe\Big(k\frac{y-\overline x}{q}\Big)
                &=\frac{\phi(q)}{q}\cdot
                        \frac{|\mathfrak{I}(\Omega)|}{2}\cdot |\I|+O\left(q^{1/2+\varepsilon}\right).
        \end{split}
  \end{equation}  
The required result follows by \eqref{eqNFeI1} and \eqref{eqNFeI2},
and the fact that
  \begin{equation*}
        \begin{split}
        \Neo(\Omega)=\sum_{\substack{1\le q\le R\\ q\text{ odd}}}
                        \frac{\phi(q)}{q}\cdot\frac{|\mathfrak{I}(\Omega)|}{2}+O(R)\,.
        \end{split}
  \end{equation*}  

\end{proof}

Now we have all the tools needed to complete the proof of Theorem~\ref{TheoremI}.
Since $\FQe=Q^2/\pi^2+O(Q\log Q)$, by Lemma~\ref{LemmaFQeI} we find
that $\#\FQeI=|\I|\cdot\#\FQe+O(Q^{3/2+\varepsilon})$.
On the other hand, using notations from Section~\ref{secProofT}, we find that
$\B^\I_Q(r)=\Neo^\I(\Omega_Q(r))$, where
$\Omega_Q(r)=\Omega_Q(x_0,y_0,\eta)(r)$ is given by 
  \begin{equation*}
     \Omega_Q(r)=\left\{(x,y)\in \RR^2\colon\ 
        \begin{array}{l}  
              1\le x,y\le Q,\  x+y>Q,\  \\ \displaystyle
         \kk(x,y)\in\A(r),\     (x,x_r^\LL(x,y))\in Q\cdot\Box 
        \end{array} 
        \right\}.
  \end{equation*}
The set $\Omega_Q(r)$ is in general not convex, but it is a finite
union of boundedly many convex sets, as $Q\to\infty$, the number of
these convex sets depending on the given point
$(x_0,y_0)\in [0,1]\times [0,1]$. Then, by
Lemma~\ref{LemmaNeoI}, it follows that 
$$\#\B_Q^\I(r)= |\I|\cdot\#\B_Q(r)+O_{(x_0,y_0)}(Q^{3/2+\varepsilon}),$$
concluding the proof of Theorem~\ref{TheoremI}.


\end{document}